\documentclass[12pt]{amsart}
\usepackage{amssymb,latexsym}
\newdimen\AAdi%
\newbox\AAbo%
\def\AAk#1#2{\s_etbox\AAbo=\hbox{#2}\AAdi=\wd\AAbo\kern#1\AAdi{}}%
\def\AAr#1#2#3{\s_etbox\AAbo=\hbox{#2}\AAdi=\ht\AAbo\raise#1\AAdi\hbox{#3}}%
\font\tenmsb=msbm10 at 12pt \font\sevenmsb=msbm7 at 8pt
\font\fivemsb=msbm5 at 6pt
\newfam\msbfam
\textfont\msbfam=\tenmsb \scriptfont\msbfam=\sevenmsb
\scriptscriptfont\msbfam=\fivemsb
\def\Bbb#1{{\tenmsb\fam\msbfam#1}}
\begin{document}
\newcommand{\D}{D \hskip -2.8mm \slash}
\newcommand{\pp}{\overline{\psi}}
\newcommand{\e}{\overline{\varepsilon}}
\newcommand{\hh}{\sqrt{h}d^2x}
\newcommand{\ii}{\frac{1}{2}\int_M}
\renewcommand{\theequation}{\thesection.\arabic{equation}}
\newcommand{\wb}{\widetilde{\nabla}_{e_\beta}}
\newcommand{\wa}{\widetilde{\nabla}_{e_\alpha}}
%\numberwithin{equation}{section}e_
\newtheorem{thm}{Theorem}
\newtheorem{lem}{Lemma}
\newtheorem{cor}{Corollary}
\newtheorem{rem}{Remark}
\newtheorem{pro}{Proposition}
\newtheorem{defi}{Definition}
\newcommand{\noi}{\noindent}
\newcommand{\dis}{\displaystyle}
\newcommand{\mint}{-\!\!\!\!\!\!\int}
\newcommand{\ba}{\begin{array}}
\newcommand{\ea}{\end{array}}
\def \tf{\tilde{f}}
\def\cqfd{%
\mbox{ }%
\nolinebreak%
\hfill%
\rule{2mm} {2mm}%
\medbreak%
\par%
}
\def \pr {\noindent {\it Proof.} }
\def \rmk {\noindent {\it Remark} }
\def \esp {\hspace{4mm}}
\def \dsp {\hspace{2mm}}
\def \ssp {\hspace{1mm}}
\def\n{\nabla}
\def\RR{\Bbb R}\def\R{\Bbb R}
\def\C{\Bbb C}
\def\B{\Bbb B}
\def\N{\Bbb N}
\def\Q{\Bbb Q}
\def\Z{\Bbb Z}
\def\EE{\Bbb E}
\def\H{\Bbb H}
\def\SS{\Bbb S}\def\S{\Bbb S}
\def \c {{\bf C}}
\def \Z {{\bf Z}}
\def \Q {{\bf Q}}
\def \a {\alpha}
\def \b {\beta}
\def \d {\delta}
\def \e {\epsilon}
\def \G {\Gamma}
\def \g {\gamma}
\def \l {\lambda}
\def \L {\Lambda}
\def \O {\Omega}
\def \om {\omega}
\def \o{\omega}
\def \s {\sigma}
\def \t {\theta}
\def \z {\zeta}
\def \vp {\varphi}
\def \vt {\vartheta}
\def \ve {\varepsilon}
\def \i {\infty}
\def \ds {\displaystyle}
\def \oo {\overline{\Omega}}
\def \ov {\overline}
\def \bd {\bigtriangledown}
\def \U {\bigcup}
\def \un {\underline}
\def \h {\hspace{.5cm}}
\def \hs {\hspace{2.5cm}}
\def \v {\vspace{.5cm}}
\def \mi {M_{i}}
\def \ra {\longrightarrow}
\def \Ra {\Longrightarrow}
\def \rw {\rightarrow}
\def \bs {\backslash}
\def \rn {{\bf R}^n}
\def \h* {\hspace*{1cm}}
\def\la{\langle}
\def\ra{\rangle}

\def \vs{\vspace*{0.1cm}}
\def\hs{\hspace*{0.6cm}}
\def \ds{\displaystyle}
\def\cal{\mathcal}
\title[Regularity Theorems for Dirac-Harmonic Maps]{ Regularity Theorems and Energy Identities for Dirac-Harmonic Maps}
\author{Qun Chen, J\"urgen Jost, Jiayu Li, Guofang Wang}
\thanks{The research of QC and JYL was partially supported by  NSFC}
\address{School of
Mathematics and Statistics\\ Central China Normal University\\
Wuhan 430079, China } \email{qunchen@mail.ccnu.edu.cn}
\address{Max Planck Institute for Mathematics in the Sciences \\Inselstr. 22-26\\D-04103 Leipzig, Germany
} \email{jjost@mis.mpg.de}
\address{Partner Group of Max Planck Institute for Mathematics in the Sciences\\
Institute of Mathematics\\
Chinese Academy of Sciences\\
Beijing 100080, P. R. of China} \email{lijia@mail.amss.ac.cn}
\address{Max Planck Institute for Mathematics in the Sciences\\Inselstr. 22-26\\D-04103 Leipzig, Germany
} \email{gwang@mis.mpg.de} 
\begin{abstract}
We study Dirac-harmonic maps from a Riemann surface to a sphere
$\S^n$. We
   show that a weakly Dirac-harmonic map is in fact smooth,
   and prove that the energy identity
   holds during the blow-up process.
\end{abstract}
%\date{Mar. 11, 2004}
\maketitle
\section{Introduction}
Let $M$ be a compact spin Riemann surface, $\Sigma M$ the spinor
bundle over $M$ and $N$ a compact Riemannian manifold. Let $\phi$ be
a map from $M$ to $N$, $\psi$ a section of the bundle $\Sigma M
\otimes\phi^{-1}TN$. Let $\widetilde{\nabla}$ be the connection
induced from those on $\Sigma M$ and $\phi^{-1}TN$. The Dirac
operator $\D$ along the map $\phi$ is defined by $\D
\psi:=e_\a\cdot\widetilde{\nabla}_{e_\a}\psi,$ where $e_1,e_2$ is
an orthonormal basis on $M$. We consider the functional
$$L(\phi,\psi):=\int_M[|d\phi|^2+\la\psi,\D\psi\ra_{\Sigma
M\otimes TN} ].$$ The critical points $(\phi,\psi)$ are called
Dirac-harmonic maps from $M$ to $N$ (these maps were first introduced
in our companion paper 
\cite{CJLW1} were also further background and motivation are
provided). When $\psi$ vanishes, we obtain the 
standard energy functional whose minimizers $\phi$ are harmonic
maps. In other words, here we are generalizing that setting by
coupling the map with a spinor field with values in the pull-back
tangent bundle. The important point is that this generalization
preserves a fundamental property of the energy functional on
Riemann surfaces, namely its conformal invariance. In fact, our
functional is nothing but the action functional for the non-linear
supersymmetric sigma model from quantum field theory, with the
only difference that here all fields are real valued instead of
having Grassmann coefficients. This brings us back into the
framework of the calculus of variations.
\par
Since the construction is geometrically quite natural, one should
expect that this class of  maps can yield new geometric invariants
of $N$. Before one can address that issue, however, one needs to
do the basic analytic work. As a first step, we should derive a
compactness theorem. To begin this program, we consider in this
paper the  case that the target is a sphere $\S^n$, that is, in
the terminology of quantum field theory, we consider the O$(n+1)$
sigma model. Suppose that $(\phi_k,\psi_k)$ is a sequence of
Dirac-harmonic maps from $M$ to $\S^n$ with uniformly bounded
energy $E(\phi_k,\psi_k)=\int_M(|d\phi_k|^2+|\psi_k|^4)$, then
there is a subsequence which we also denote by $(\phi_k,\psi_k)$
such that $\phi_k\to \phi$ weakly in $W^{1,2}$ and $\psi_k\to\psi$
weakly in $L^4$, and outside a finite set of points
$S=\{p_1,p_2,\cdots,p_I\}$ which we call the blow-up set, the
convergence is strong on compact sets. So $(\phi,\psi)$ is smooth
in $M\setminus S$, and it is a weakly Dirac-harmonic map. We show
first in this paper that any weakly Dirac-harmonic map is smooth,
and so also the present limit is smooth. At every blow-up point
$p_i$, by Sacks-Uhlenbeck's blow-up, one gets a finite
number of Dirac-harmonic spheres $(\sigma^l_i,\xi^l_i)$. 
The regularity of weakly harmonic maps  was proved by Helein \cite{H2} and 
\cite{H1}.
Another
main purpose of this paper is to show the so-called energy
identity: $${\rm lim}_{k\to \infty}
E(\phi_k)=E(\phi)+\Sigma_i\Sigma_l E(\sigma^l_i);\qquad {\rm
lim}_{k\to \infty} E(\psi_k)=E(\psi)+\Sigma_i\Sigma_l
E(\xi^l_i).$$ The energy identity for a min-max sequence for the
energy was proved by Jost \cite{J2}, for Palais-Smale sequences with
uniformly $L^{2}$-bounded tension fields by Ding-Tian \cite{DT}. 
For related results see  \cite{PW}, \cite{Y}, \cite{P}, \cite{QT}  and
\cite{Lamm}.
%for harmonic maps it was
%proved by Parker-Wolfson and  \cite{PW} Ye \cite{Y}.
\par

\section{Regularity theorems for Dirac-harmonic maps}
Let $(M,h_{\a\b})$ be a compact two-dimensional Riemannian manifold
with a fixed 
spin structure, $\Sigma M$ the spinor bundle. For any $X\in \Gamma
(TM),$ $\xi\in\Gamma (\Sigma M),$ denote by $X\cdot\xi$ the
Clifford multiplication, which satisfies the following
skew-adjointness relation: $$\la X\cdot\xi, \eta\ra_{\Sigma
M}=-\la \xi, X\cdot\eta\ra_{\Sigma M}$$ for any $X\in \Gamma
(TM),$ $\xi,\eta \in\Gamma (\Sigma M),$ where $\la
\cdot,\cdot\ra_{\Sigma M}$ denotes the metric on $\Sigma M$ induced by
the Riemannian metric $h_{\a\b})$.
Choosing a local orthonormal basis $\{e_\alpha, \alpha=1,2\}$ on
$M$, the usual Dirac operator is defined as:
 $
\partial \hskip -2.2mm \slash
:=e_{\a}\cdot \n_{e_\a},$ where $\n$ stands for the spin
connection on $\Sigma M$ (here and in the sequel, we use the
Einstein summation convention). (A good reference for the spin
geometry tools used in this paper is \cite{LM}.)\\  
Let $\phi$ be a
smooth map from $M$ to another compact Riemannian manifold $(N,g)$
of dimension $n\ge 2$. Let $\phi^{-1} TN$ be the pull-back bundle of
$TN$ by $\phi$ and consider the twisted bundle $\Sigma M\otimes
\phi^{-1} TN$. On $\Sigma M\otimes \phi^{-1} TN$ there is a  metric
$\la\cdot,\cdot\ra_{\Sigma M\otimes TN}$ induced from the metrics
on $\Sigma M$ and $\phi^{-1}TN$.  Also we have a natural connection
$\widetilde{\n}$ on $\Sigma M\otimes \phi^{-1} TN$ induced from those
on $\Sigma M$ and $\phi^{-1} TN$. In local coordinates, the section
$\psi$ of $\Sigma M\otimes \phi^{-1} TN$ is written as
\[\psi= \psi^j\otimes\partial_{y^j}(\phi),\]
where each $\psi^j$ is a usual spinor on $M$ and $\{\partial_{
y^j}\}$ is the natural local basis on $N$.  $\widetilde{\n}$
becomes
\[\widetilde{\nabla} \psi= \n
\psi^i\otimes\partial_{y^i}(\phi)+ (\Gamma^i_{jk} \n \phi^j)
\psi^k\otimes\partial_{y^i}(\phi)\] where the $\Gamma^i_{jk}$ are the
Christoffel symbols of the Levi-Civita connection of $N$. \\
We define
the {\it Dirac 
operator along the map $\phi$} as
 \begin{eqnarray*}
\D\psi&:=& e_\a\cdot \widetilde{\nabla}_{e_\a} \psi\nonumber
\\ &=&
\partial \hskip -2.2mm \slash
 \psi^i\otimes \partial_{y^i}(\phi)
+ (\Gamma^i_{jk} \n_{e_\a} \phi^j) (e_\a\cdot
\psi^k)\otimes\partial_{y^i}(\phi).
\end{eqnarray*}
It is easy to verify that $\D$ is formally self-adjoint, i.e., $$
 \int_M\la\psi, \D\xi\ra_{\Sigma
M\otimes TN}=\int_M\la\D\psi,\xi\ra_{\Sigma M\otimes TN},
 $$
for all $\psi, \xi \in C^\infty (\Sigma M\otimes \phi^{-1}TN)$, the
space of smooth sections of $\Sigma M\otimes \phi^{-1}TN$, where
$\la\psi,\xi\ra_{\Sigma M\otimes
TN}:=g_{ij}(\phi)\la\psi^i,\xi^j\ra_{\Sigma M}$, for $\psi,\xi\in
\Gamma (\Sigma M\otimes \phi^{-1}TN)$. \\
Set
\begin{equation}
\label{aux1}
{\cal X}(M,N):=\{ (\phi,\psi)\,|\, \phi \in C^\infty(M,N) \hbox{ and }
\psi \in C^\infty (\Sigma M\otimes \phi^{-1}TN)\}.
\end{equation}
We consider
the following functional
\begin{eqnarray*}
 L(\phi,\psi)&:=&\int_M (|d\phi|^2+\la\psi,\D\psi\ra_{\Sigma
M\otimes TN}) \nonumber \\ &=&\int_M
\{g_{ij}(\phi)h^{\alpha\beta}\frac{\partial \phi^i}{\partial
x^\alpha}\frac{\partial \phi^j}{\partial
x^\beta}+g_{ij}(\phi)\la\psi^i,\D\psi^j\ra_{\Sigma M}\}.
\end{eqnarray*}
By a direct computation, we obtain the Euler-Lagrange equations of
$L$:
\begin{equation}
\label{b1} \D\psi^i=\partial \hskip -2.2mm \slash
\psi^i+\Gamma^i_{jk}(\phi)\partial_\alpha\phi^j(e_\alpha\cdot\psi^k)=0,\quad
i=1,2,\cdots,n,
\end{equation}
\begin{equation}
\label{b2}
\tau^m(\phi)-\frac{1}{2}R^m\hskip-1mm_{lij}(\phi)\la\psi^i,\nabla\phi^l\cdot\psi^j\ra_{\Sigma
M}=0, \quad m=1,2,\cdots,n,
\end{equation}
where $\tau(\phi)$ is the tension field of the map $\phi$,
$\n\phi^l\cdot\psi^j$ denotes the Clifford multiplication of the
vector field $\n\phi^l$ with the spinor $\psi^j$, and
$R^m\hskip-1mm_{lij}$ stands for a component of the curvature
tensor of the target manifold $N$. % (In this paper, we use the same
%curvature convention as in \cite{CE}.)
Denote $$\cal R(\phi,\psi):=\frac 1 2
R^m\hskip-1mm_{lij}\la\psi^i,\nabla \phi^l\cdot\psi^j\ra\partial
_{y_m}.$$
 We write equations (\ref{b1}) (\ref{b2}) in global form as 
\begin{equation}
\label{b1n} \begin{array}{rcl}
\D\psi &=&\vs 0\\
\tau(\phi)&=&{\cal R}(\phi,\psi).\end{array}
\end{equation}

 Solutions $(\phi,\psi)$ of 
(\ref{b1}) and (\ref{b2}) are called {\it Dirac-harmonic maps from
$M$ to $N$}. When $\psi=0$, a solution $(\phi,0)$ is just a harmonic map.
Harmonic maps have been extensively studied. See, for instance, two reports
of Eells-Lemaire \cite{EL}.
When $\phi$ is a constant map, each component
of $\psi$ is a usual harmonic spinor.  Harmonic spinors also have
been well understood, see for instance \cite{Hi}, \cite{LM}, \cite{BS} and \cite
{B}.
%In \cite{Hi}, Hitchin
%showed that the dimension of harmonic spinors is conformally
%invariant, B\"{a}r and Schmutz determined in \cite{BS} the
%dimension of harmonic spinors in the hyperelliptic case, Maier discussed
%generic metrics in \cite{St}, these works enlarged the knowledge about harmonic
%spinors on Riemann surfaces. Our
Dirac-harmonic maps thus are a generalization and combination of
harmonic maps and  
harmonic spinors. Non-trivial examples are given in \cite{CJLW1}.
\par
Let $(N', g')$ be another Riemannian manifold and $f: N\to N'$ a
smooth map. %In the following we simply write $f_*$ to denote the
%map $id\otimes f_*: \Sigma M\otimes\phi^{-1}TN \to \Sigma
For a map $\phi: M\to N$, we have a map $\phi'=\phi\circ f$ from
$M$ to $N'$. The map $f$ naturally induces a map from
$M\otimes\phi^{-1}TN \to M\otimes\phi'^{-1}TN,$which is denoted by $f_*$.
Hence for and $(\phi, \psi)\in {\cal X}$ we get 
$(\phi', f_*\psi)\in {\cal X}(M,N')$.  $\psi ':=f_*\psi$ is a spinor field along the
 map $ \phi'$.

\par
 Let $A$ be
the second fundamental form of $f$, i.e., $A(X,Y)=( \n_X df)(Y)$
for any $X,Y\in \Gamma(TN)$. 
It is well-known that the tension fields of $\phi$ and $\phi '$ satisfy the
following relation
\begin{equation}
\label {b3}
  \tau' ( \phi ')=A(d\phi(e_\a),
d\phi(e_\a))+ df(\tau(\phi)).
\end{equation}
%For any $(\phi, \psi)\in {\cal X}$ we
%set
%\[\phi ' =f\circ\phi \quad \hbox{ and }  \quad  \psi '= f_*\psi.\]
 One can also check that the Dirac operators $\D$ and
$\D\hskip1mm '$ corresponding to $\phi$ and $\phi '$ respectively
are related by 
\begin{equation}
\label {b4}
 \D \hskip1mm ' \psi '= f_*(\D\psi)+\cal A(d\phi(e_\a), e_\a\cdot
\psi ),
\end{equation}
where
\begin{eqnarray*}
\cal A(d\phi(e_\a), e_\a\cdot \psi )&:=&\phi^i_\a\hskip0.1mm
e_\a\cdot\psi^j\otimes A(\partial_{y^i},\partial_{y^j})
\\&=&(\n\phi^i\cdot\psi^j)\otimes A(\partial_{y^i},\partial_{y^j}).
\end{eqnarray*}

When $f: N\to N'$ is an isometric immersion, then $A(\cdot,\cdot)$
is the second fundamental form of the submanifold $N$ in $N'$. We have
\[\n'_X\xi=-P(\xi;X)+\n_X^{\perp}\xi, \qquad \quad \n'_XY=\n_XY+A(X,Y)  \]
$\forall X,Y\in \Gamma(TN),$ $\xi\in\Gamma (T^{\perp}N)$, where
$T\perp N$ is the normal bundle, $\nabla$ and $\nabla'$ are
covariant derivatives and
$P(\cdot\hskip1mm ; \cdot)$ denotes the shape operator. 
In this case, for simplicity of notation, we identify
$\phi$ with $\phi'$ and $\psi$ with $\psi'$.
Using the equation of Gauss, we have
\begin{eqnarray*}
&&R^m\hskip-1mm_{lij}\la\psi^i,\nabla \phi^l\cdot\psi^j\ra_{\Sigma
M}\\ &=&g^{mk}[\la A(\partial_{y^k},\partial_{y^i}),
A(\partial_{y^l},\partial_{y^j})\ra_{TN'}- \la
A(\partial_{y^k},\partial_{y^j}),
A(\partial_{y^l},\partial_{y^i})\ra_{TN'}\\
&&+R'_{klij}]\la\psi^i,\nabla \phi^l\cdot\psi^j\ra_{\Sigma M}\\
&=&2g^{mk}\la A(\partial_{y^k},\partial_{y^i}),
A(\partial_{y^l},\partial_{y^j})\ra_{TN'}\la\psi^i,e_\alpha\cdot\psi^j\ra_{\Sigma
M}\phi^l_\alpha
\\ &&+g^{mk}R'_{klij}\la\psi^i,\nabla \phi^l\cdot\psi^j\ra_{\Sigma M} \\
&=&2g^{mk}\la
P(A(\partial_{y^l},\partial_{y^j});\partial_{y^i}),\partial_{y_k}\ra_{TN}\la\psi^i,e_\a\cdot\psi^j\ra_{\Sigma
M}\phi^l_\a +R'\hskip0.02mm^m\hskip-1mm_{lij}\la\psi^i,\nabla
\phi^l\cdot\psi^j\ra_{\Sigma M},
\end{eqnarray*}
where in the last step we used the following relation between the
shape operator $P(\cdot\hskip1mm;\cdot)$ and the second
fundamental form $A(\cdot,\cdot)$: $$ \la P(\xi;X),Y\ra_{TN}=\la
A(X,Y),\xi\ra_{TN'}$$ for any $X,Y\in \Gamma (TN),$ $\xi\in \Gamma
(T^{\perp}N)$.
Set $$P(\cal
A(d\phi(e_\alpha),e_\alpha\cdot\psi);\psi):=P(A(\partial_{y^l},\partial_{y^j});\partial_{y^i})\la\psi^i,e_\a\cdot\psi^j\ra_{\Sigma
M}\phi^l_\a. $$ From the above calculation, we have
\begin{equation}
\label {b5} \cal R(\phi,\psi)=P(\cal
A(d\phi(e_\alpha),e_\alpha\cdot\psi);\psi)+\cal R'(\phi,\psi).
\end{equation}
Therefore, using (\ref{b4}) and (\ref{b3}) and identifying
$\psi$ with $ \psi '$ and $\phi$ with $ \phi '$, we can rewrite
(\ref{b1}) and (\ref{b2})  as follows:
\begin{eqnarray}
\label{b6}
 \D \hskip1mm ' \psi & =& \cal A(d\phi(e_\a), e_\a\cdot \psi),\\
 \label{b7} \tau' ( \phi ) &=& A(d\phi(e_\a), d\phi(e_\a))+P(\cal
A(d\phi(e_\alpha),e_\alpha\cdot\psi);\psi)+ \cal R'(\phi,\psi).
\end{eqnarray}
In particular, by the Nash-Moser embedding theorem, we embed $N$
into the Euclidean space $N'=\RR^K$, and have $\D\hskip1mm '=\partial
\hskip -2.2mm \slash$ and $\tau'=-\Delta$, where
$\Delta$ is the (negative) Laplacian.
Therefore, we have  
\begin{eqnarray}
\label{b8}
\partial \hskip -2.2mm \slash
 \psi &=& \cal A(d\phi(e_\a), e_\a\cdot \psi ),\\
%\end{equation}
%\begin{equation}
\label{b9} -\Delta \phi &=& A(d\phi, d\phi)+P(\cal
A(d\phi(e_\alpha),e_\alpha\cdot\psi);\psi),
\end{eqnarray}
where $\phi:M\to \R^K$ with 
\begin{equation}\label{add1}
\phi(x)\in N \end{equation}
for any $x\in M$
 and
$\psi=(\psi^1,\psi^2,\cdots, \psi^K)$
%:=\psi^i\otimes E_i$ (we
%always use $\{E_i|i=1,2,\cdots,K\}$ to denote the standard
%orthonormal basis of $\R^K$) 
with the property that $\psi(x)$ is
along the map $\phi$, namely,
\begin{equation}\label{eq100}
\sum_{i=1}^K v_i\psi^i(x)=0, \quad \hbox{ for any normal vector }
v=(v_1,\cdots, v_K) \hbox{ at } \phi(x).\end{equation} Here
$\psi^i\in \Gamma (\Sigma M)$. For any vector
$v=(v_1,v_2,\cdots,v_K) \in \R^K$, abusing the notation a little bit,
we write $\la v,
\psi\ra=v_i\psi^i\in \Gamma (\Sigma M)$. And we also write 
for $\xi\in \Gamma(\Sigma M)$
\[\la\psi, \xi\ra :=(\la \psi^1, \xi\ra_{\Sigma M}, \la \psi^2, \xi\ra_{\Sigma M}, \cdots,
\la \psi^K, \xi\ra_{\Sigma M})\in \R^K,\]
if there is no confusion.

%Now using (\ref{b8}) and (\ref{b9}), we can define weakly Dirac-harmonic maps. 
Set
\[{\cal X}^{1,2}_{1,4/3}(M,N):=
\{(\phi,\psi)\in W^{1,2}\times W^{1,4/3} \hbox{ with } (\ref{add1}) \hbox{ and }
(\ref{eq100}) \hbox{ a.e.}\}.\]
For simplicity of notation, we denote ${\cal X}^{1,2}_{1,4/3}(M,N)$
by ${\cal X}(M,N)$ (thus, we are changing the convention of (\ref{aux1})).
It is clear that the functional $L(\phi,\psi)$ 
is well-defined for $(\phi,\psi)\in {\cal X}(M,N)$. 

\vskip12pt 
\noindent {\bf Definition.} {A critical point $(\phi,\psi)\in 
{\cal X}(M,N)$ of the functional $L$ in ${\cal X}(M,N)$ is called
%Suppose
%that $(\phi,\psi)\in {\cal X}(M,N)$.
%We call $(\phi,\psi)$
 {\it a weakly Dirac-harmonic map}
from $M$ to $N$. Equivalently, $(\phi,\psi)\in 
{\cal X}(M,N)$  is  a weakly Dirac-harmonic map
from $M$ to $N$ if and only if $(\phi,\psi)$ satisfies
\begin{eqnarray}
\label{b10} \int_M\{\la\nabla \phi,\nabla \eta\ra-\la
A(d\phi,d\phi)+P(\cal
A(d\phi(e_\alpha),e_\alpha\cdot\psi);\psi),\eta\ra\}&=& 0,\\
\label{b11} \int_M\{\la \psi,\partial \hskip -2.2mm \slash \xi
\ra-\la \cal A(d\phi(e_\a), e_\a\cdot \psi ),\xi\ra\}&=&0,
\end{eqnarray}
for all $\eta \in C^\infty(M,\RR^K)$ and $\xi \in C^\infty(\Sigma
M\otimes \RR^K)$.
 }
\vskip12pt   One of our purposes of this paper 
is to study the regularity of weakly
Dirac-harmonic maps.
Our main observation is  that when the target $N$ is the standard
sphere $\S^n$, a weakly Dirac-harmonic map has a special structure
like a weakly harmonic map. For  weakly harmonic maps, see \cite{H1}.
 
 \vskip12pt
 \noindent
{\bf Proposition 2.1.}  {\it Let $M$ be a Riemann surface with a fixed
spin structure and $(\phi,\psi)\in {\cal X}(N,M)$  a weakly Dirac-harmonic
map from $M$ to $\S^n$. Let $D$ be a simply connected domain of $M$.
Then there exists $M=(M^{ij})\in W^2(D, \R^{n\times n})$ such that
\begin{equation}\label{eqn}
-\Delta \phi=\frac{\partial M}{\partial x}\frac{\partial \phi}{\partial y}-
\frac{\partial M}{\partial y}\frac{\partial \phi}{\partial x}.\end{equation}
 }
\vskip12pt \noindent \pr   For $N=\S^n\subset \RR^{n+1}$,   the
equations (\ref{b8}) and (\ref{b9}) can be respectively written as
follows:
\begin{eqnarray}
\label{b12}
\partial \hskip -2.2mm \slash
 \psi^m & =& -\sum(\n\phi^i\cdot\psi^i)
 %\la d\phi(e_\a), e_\a\cdot \psi \ra
 \otimes\phi^m,\\
%\end{equation}
%\begin{equation}
\label{b13} -\Delta \phi^m &=&|d\phi|^2\phi+\la \psi^m\otimes
d\phi(e_\alpha),e_\alpha\cdot\psi\ra_{\Sigma M\otimes \R^K},
\end{eqnarray}
for $m=1,2,\cdots, n+1$.
%Given any point $x_0\in M$, we prove that $\phi$ is continuous at
%$x_0$, since this is a local property, and the functional is
%conformally invariant under conformal diffeomorphisms of this
%domain, we may choose a unit disk $D$ in $\RR^2$ with standard
Set $\phi_\a:=d\phi(e_\a)$.
\par From (\ref{b13}), we have for $m=1,2,\cdots,n+1$ that
\begin{equation}\label{b14}
\begin{array}{rcl}
 \Delta \phi^m&=&\vs\ds -|d \phi|^2\phi^m+\la \psi^m\otimes
d\phi(e_\alpha),e_\alpha\cdot\psi\ra_{\Sigma M\otimes
\R^K} \\
&=&\vs\ds -(\phi^i_\a\phi^m-\phi^i\phi^m_\a)\phi^i_\a+\phi^i_\a\la
e_\a\cdot\psi^i,\psi^m\ra_{\Sigma M}\quad (\hbox{as }
0=\partial_\alpha 1=2\phi^i_\alpha\phi^i)
\\
&=&\ds\vs [\la\partial_x\cdot\psi^i,\psi^m\ra_{\Sigma
M}-(\phi^i_x\phi^m-\phi^i\phi^m_x)]\phi^i_x \\
&&+[\la\partial_y\cdot\psi^i,\psi^m\ra_{\Sigma
M}-(\phi^i_y\phi^m-\phi^i\phi^m_y)]\phi^i_y
\\ 
&:=& \ds A^{mi}\phi^i_x+B^{mi}\phi^i_y.
\end{array}\end{equation}
We would like to show that there exists some function $M^{mi}$ on
$D$ such that
\begin{equation}
\label{b15} A^{mi}=M^{mi}_y,\quad B^{mi}=-M^{mi}_x
\end{equation}
which, by the Frobenius theorem, is equivalent to
\begin{equation}
\label{b16} A^{mi}_x+B^{mi}_y=0.
\end{equation}
Calculating directly, one derives
\begin{eqnarray}
\label{b17}
 A^{mi}_x+B^{mi}_y&=&
\la\partial_x\cdot\psi^i_x,\psi^m\ra_{\Sigma
M}+\la\partial_x\cdot\psi^i,\psi^m_x\ra_{\Sigma M}
+\la\partial_y\cdot\psi^i_y,\psi^m\ra_{\Sigma M}\nonumber \\
&&+\la\partial_y\cdot\psi^i,\psi^m_y\ra_{\Sigma M}
-(\phi^i_{xx}\phi^m-\phi^i\phi^m_{xx}+\phi^i_{yy}\phi^m-\phi^i\phi^m_{yy})\nonumber
\\ &=&\la\partial \hskip -2.2mm \slash
\psi^i,\psi^m\ra_{\Sigma M}-\la\psi^i,\partial \hskip -2.2mm
\slash\psi^m\ra_{\Sigma
M}-(\Delta\phi^i\phi^m-\Delta\phi^m\phi^i).
\end{eqnarray}
 From equation (\ref{b13}), one gets
\begin{equation}
\label{b18}
-(\Delta\phi^i\phi^m-\Delta\phi^m\phi^i)=-\phi^k_\alpha\phi^m\la
e_\alpha\cdot\psi^k,\psi^i\ra_{\Sigma M}+\phi^k_\alpha\phi^i\la
e_\alpha\cdot\psi^k,\psi^m\ra_{\Sigma M},
\end{equation}
And  from  equation (\ref{b12}), one obtains
$$\la\partial \hskip -2.2mm \slash \psi^i,\psi^m\ra_{\Sigma
M}=-\phi^k_\alpha\phi^i\la e_\alpha\cdot\psi^k,\psi^m\ra_{\Sigma
M}.$$ Hence
\begin{equation}
\label{b19} \la\partial \hskip -2.2mm \slash
\psi^i,\psi^m\ra_{\Sigma M}-\la\psi^i,\partial \hskip -2.2mm
\slash \psi^m\ra_{\Sigma M}=\phi^k_\alpha\phi^m\la
e_\alpha\cdot\psi^k,\psi^i\ra_{\Sigma M}-\phi^k_\alpha\phi^i\la
e_\alpha\cdot\psi^k,\psi^m\ra_{\Sigma M}.
\end{equation}
Putting (\ref{b18}) and (\ref{b19}) into (\ref{b17}) we have
(\ref{b16}). Hence we prove the Proposition. \hfill $\Box$

\
 
 \noindent{\bf Theorem 2.2.} {\it Let $M$ be a Riemann surface with a fixed
spin structure. Suppose that $(\phi,\psi)\in {\cal X}(N,M)$ 
is  a weakly Dirac-harmonic
map from $M$ to $\S^n$. Then $\phi\in C^0$. And hence $(\phi,\psi)$ is
smooth.}

\

\noindent{\it Proof.} From Proposition 2.1 and 
 Wente's well-known lemma (\cite{W}), we know that $\phi^m$ is
continuous, $m=1,2,\cdots,$ $n+1$, namely, $\phi\in
C^0(M,\S^n)$. By Theorem 2.3 below, we have that $\phi$ and $\psi$ are
smooth.\hfill $\Box$

\par

  \vskip12pt \noindent {\bf Theorem 2.3.} {\it
Let $(\phi,\psi): (D,\delta_{\alpha\beta})\to (N^n,g_{ij})$ be a
weakly Dirac-harmonic map. If $\phi$ is continuous, then
$(\phi,\psi)$ is smooth. } \vskip12pt To prove this theorem, we
first establish two lemmas (Lemma 2.4 and Lemma 2.5 below which
are similar to Lemma 8.6.1 and Lemma 8.6.2 in \cite{J1}). Since
$\phi$ is in $C^0(D,N)$, we can choose local coordinates $\{y_i\}$
on $N$ such that $\Gamma^i_{jk}(\phi(0))=0$. In these coordinates,
the equations for $\phi$ and $\psi$ can be written as
\begin{equation}
\label{d1}
\Delta\phi^m=-\Gamma^m_{ij}(\phi)\phi^i_\alpha\phi^j_\alpha+\frac
1 2
R^m\hskip-1mm_{jkl}(\phi)\la\psi^k,\nabla\phi^j\cdot\psi^l\ra_{\Sigma
M},
\end{equation}
\begin{equation}
\label{d2}
\partial \hskip -2.2mm \slash
\psi^m=-\Gamma^m_{ij}(\phi)\nabla\phi^i\cdot\psi^j.
\end{equation}
This intrinsic version of our equations is well-defined since
$\phi\in C^0$. \vskip12pt \noindent {\bf Lemma 2.4.}  {\it Let
$(\phi,\psi)$ be a weak solution of (\ref{d1}) and (\ref{d2}). If
$\phi\in C^0\cap W^{1,2}(D,N)$, then for any $\varepsilon >0$,
there is a $\rho >0$ such that
\begin{equation}
\label{d3} \int_{D(x_1,\rho)}|\nabla\phi|^2\eta^2(x)\leq
\varepsilon\int_{D(x_1,\rho)}|\nabla\eta|^2+C\varepsilon(\int_{D(x_1,\rho)}
|\psi|^4\eta^4)^{\frac 1 2},
\end{equation}
where  $D(x_1,\rho)\subset D$, $\eta\in
W_0^{1,2}(D(x_1,\rho),\RR)$, $C$ is a positive constant
independent of $\varepsilon,\rho, \phi$ and $\psi$.} \vskip12pt
\noindent {\it Proof.} Denote
$$G^m(x,\phi,d\phi):=\Gamma^m_{ij}(\phi)\phi^i_\alpha\phi^j_\alpha-\frac
1 2
R^m\hskip-1mm_{jkl}(\phi)\la\psi^k,\nabla\phi^j\cdot\psi^l\ra_{\Sigma
M},\quad G=(G^1,G^2,\cdots,G^n),$$ then $$|G_x|\leq
C(|d\phi|^3+|\nabla\psi||\psi||d\phi|),$$ $$|G_\phi|\leq
C(|d\phi|^2+|\psi|^2|d\phi|),$$ $$|G_{d\phi}|\leq
C(|d\phi|+|\psi|^2).$$ The weak form of equation (\ref{d1}) is:
\begin{equation}
\label{d4}
\int_D\nabla_\alpha\phi^i\nabla_\alpha\zeta^i=\int_DG^i(x,\phi,d\phi)\zeta^i,
\end{equation}
for any $\zeta\in W^{1,2}\cap L^\infty(D,\RR^K).$ Now we choose
$\zeta(x)=(\phi(x)-\phi(x_1))\eta^2(x)$, then
\begin{equation}
\label{d5}
\int_{D(x_1,\rho)}\nabla_\alpha\phi^i\nabla_\alpha\zeta^i=\int_{D(x_1,\rho)}|d\phi|^2\eta^2+2
\int_{D(x_1,\rho)}\eta(\phi^i(x)-\phi^i(x_1))\nabla_\alpha\phi\nabla_\alpha\eta.
\end{equation}
We have
\begin{eqnarray}
\label{d6}
\int_{D(x_1,\rho)}G^i(x,\phi,d\phi)\zeta^i&=&\int_{D(x_1,\rho)}[\Gamma^i_{jk}(\phi)\phi^j_\alpha
\phi^k_\alpha(\phi^i(x)-\phi^i(x_1))\eta^2(x)]\nonumber \\
  &&-\frac 1 2
\int_{D(x_1,\rho)}[R^i_{jkl}(\phi)\la\psi^k,\nabla\phi^j\cdot\psi^l\ra_{\Sigma
M}(\phi^i(x)-\phi^i(x_1))\eta^2(x) ]\nonumber \\ &\leq& C_N
\varepsilon_1 {\rm
Sup}_{D(x_1,\rho)}|\phi(x)-\phi(x_1)|\int_{D(x_1,\rho)}|d\phi|^2\eta^2\nonumber
\\ &&+C_N {\rm Sup}_{D(x_1,\rho)}
|\phi(x)-\phi(x_1)|\int_{D(x_1,\rho)}|d\phi||\psi|^2\eta^2\nonumber
\\
&\leq&C_N \varepsilon_1
{\rm Sup}_{D(x_1,\rho)}|\phi(x)-\phi(x_1)|\int_{D(x_1,\rho)}|d\phi|^2\eta^2\nonumber \\
&&+C_N {\rm Sup}_{D(x_1,\rho)}
|\phi(x)-\phi(x_1)|(\int_{D(x_1,\rho)}|d\phi|^2)^{\frac 1
2}\nonumber \\
&&\times (\int_{D(x_1,\rho)}|\psi|^4\eta^4)^{\frac 1 2},
\end{eqnarray}
where $\varepsilon_1>0$ is a given small number. On the other
hand,
\begin{eqnarray}
\label{d7} 2
\int_{D(x_1,\rho)}\eta(\phi^i(x)-\phi^i(x_1))\nabla_\alpha\phi\nabla_\alpha\eta
&\leq&C_N{\rm
Sup}_{D(x_1,\rho)}|\phi(x)-\phi(x_1)|\int_{D(x_1,\rho)}|d\phi||\nabla\eta|\eta \nonumber \\
&\leq& \frac 1 2 \int_{D(x_1,\rho)}|d\phi|^2\eta^2+8 {\rm
Sup}_{D(x_1,\rho)}|\phi(x)-\phi(x_1)|^2\nonumber
\\
&&\times\int_{D(x_1,\rho)}|\nabla\eta|^2.
\end{eqnarray}
Substituting (\ref{d5}), (\ref{d6}) and (\ref{d7}) into
(\ref{d4}), and choosing  $\rho$ small enough then yields
(\ref{d3}).  \hfill $\Box$
 \vskip12pt
 \noindent
 {\bf Lemma 2.5.}  {\it If $\phi \in C^0\cap W^{1,4}\cap
 W^{3,2}(D(x_0,R),N)$ and $(\phi,\psi)$ is a weak solution of
 (\ref{d1}), (\ref{d2}), then for $R$ sufficiently small, we have
 \begin{equation}
 \label{d8}
 \|\nabla^2\phi\|_{L^2(D(x_0,R/2))}+\|d\phi\|^2_{L^4(D(x_0,R/2))}\leq
 C_1\|d\phi\|_{L^2(D(x_0,R))},
 \end{equation}
 where $C_1>0$ is a constant depending on $|\phi|_{C^0(D,N)}$ and $R$.
 }
\vskip12pt
 \noindent
 {\it Proof.}
 We may assume $x_0=0\in D$. Given $\varepsilon'>0$ small, since
 $\phi$ is continuous, we can choose $R$ small enough such that
 $|\phi(x)-\phi(0)|<\varepsilon'$ for all $x\in D(x_0,R)$. For
 simplicity, we denote $B:=D(x_0,R)$.
 From equation (\ref{d2}) we have
$$|\partial \hskip -2.2mm \slash \psi|\leq |\Gamma^i_{jk}(\phi(x))-\Gamma^i_{jk}(\phi(0)|
|d\phi^j||\psi^k|\leq C_N|\phi(x)-\phi(0)||d\phi||\psi|,$$ hence
\begin{equation}
\label{d9}
 |\partial \hskip -2.2mm \slash \psi|\leq
C_N\varepsilon'|d\phi||\psi|
\end{equation}
Noting that $\partial \hskip -2.2mm \slash (\psi\eta)=\eta\partial
\hskip -2.2mm \slash \psi+\nabla\eta\cdot\psi$, we have
\begin{eqnarray}
\label{d10} \|\partial \hskip -2.2mm \slash
(\psi\eta)\|_{L^{4/3}(B)}&\leq& \|\eta \partial \hskip -2.2mm
\slash \psi\|_{L^{4/3}(B)}+\|\nabla\eta\cdot\psi\|_{L^{4/3}(B)}
\nonumber \\ &\leq&
C_N\varepsilon'\||d\phi||\psi|\eta\|_{L^{4/3}(B)}+\||\nabla\eta||\psi|\|_{L^{4/3}(B)}\nonumber
\\ &\leq&
C_N\varepsilon'\|d\phi\|_{L^2(B)}\||\psi|\eta\|_{L^4(B)}+\||\nabla\eta||\psi|\|_{L^{4/3}(B)}.
\end{eqnarray}
By the elliptic estimates for the first order equation, we have
\begin{equation}
\label{d11}
\|\nabla(\psi\eta)\|_{L^{4/3}(B)}+\|\psi\eta\|_{L^4(B)}\leq
C_R\|\partial \hskip -2.2mm \slash (\psi\eta)\|_{L^{4/3}(B)}.
\end{equation}
A proof was given in Lemma 4.8 in \cite{CJLW1}.
By choosing $R$ small enough such that
$C_N\varepsilon'\|d\phi\|_{L^2(B)}<\frac 1 2,$ we obtain from
(\ref{d10}) and (\ref{d11}) that
\begin{equation}
\label{d12}
\|\nabla(\psi\eta)\|_{L^{4/3}(B)}+\|\psi\eta\|_{L^4(B)}\leq
C_R\||\nabla\eta||\psi|\|_{L^{4/3}(B)},
\end{equation}
from which we easily derive that
\begin{equation}
\label{d13}
\||\nabla\psi|\eta\|_{L^{4/3}(B)}+\|\psi\eta\|_{L^4(B)}\leq
C_R\||\nabla\eta||\psi|\|_{L^{4/3}(B)}.
\end{equation}
For any $\zeta\in W^{1,2}_0(B,\RR^K)$,
\begin{equation}
\label{d14} \int_B \nabla\phi\nabla\zeta=-\int_B\Delta\phi
\zeta=\int_B G \zeta,
\end{equation}
Choosing $\zeta=\nabla_\gamma(\xi^2\nabla_\gamma\phi)$, where
$\xi\in C^\infty\cap W_0^{1,2}(B,\RR)$ is to be determined later,
we get
\begin{eqnarray}
\label{d15}
\int_B\nabla_\gamma(\nabla_\beta\phi)\nabla_\beta(\xi^2\nabla_\gamma\phi)&=&-
\int_B\nabla_\beta\phi\nabla_\beta(\nabla_\gamma(\xi^2\nabla_\gamma\phi))\nonumber
\\
&=&-\int_B G\nabla_\gamma(\xi^2\nabla_\gamma\phi)\nonumber
\\
&=&\int_B\nabla_\gamma G\hskip1mm\nabla_\gamma \phi\hskip1mm
\xi^2.
\end{eqnarray}
Note that
\begin{eqnarray}
\label{d16}
\nabla_\gamma(\nabla_\beta\phi)\nabla_\beta(\xi^2\nabla_\gamma\phi)
&=&|\nabla_\gamma\nabla_\beta\phi|^2\xi^2+(\nabla_\gamma\nabla_\beta\phi\nabla_\gamma\phi)
\nabla_\beta\xi^2\nonumber \\
&\geq&|\nabla^2\phi|^2\xi^2-2|\nabla^2\phi||d\phi||\xi\nabla\xi|,
\end{eqnarray}
and
\begin{eqnarray}
\label{d17} |\nabla_\gamma G\nabla_\gamma
\phi|&\leq&C_N(|d\phi|^4+|\psi||\nabla\psi||d\phi|^2+|d\phi|^3|\psi|^2\nonumber
\\
&& +|\nabla^2\phi||d\phi|^2+|\nabla^2\phi||d\phi||\psi|^2).
\end{eqnarray}
Substituting (\ref{d16}) and (\ref{d17}) into (\ref{d15}) yields
\begin{eqnarray}
\label{d18}
\int_B|\nabla^2\phi|^2\xi^2&\leq&C_N\int_B|\nabla^2\phi||d\phi||\xi\nabla\xi|+C_N
\int_B|\nabla^2\phi||d\phi|^2\xi^2\nonumber \\
&&+C_N\int_B|d\phi|^4\xi^2+C_N\int_B|\psi||\nabla\psi||d\phi|^2\xi^2\nonumber
\\
&&+C_N\int_B|d\phi|^3|\psi|^2\xi^2+C_N\int_B|\nabla^2\phi||d\phi||\psi|^2\xi^2\nonumber
\\
&:=&I+II+III+IV+V+VI.
\end{eqnarray}
For $\varepsilon_1>0$ small, we have
\begin{equation}
\label{d19} I\leq C_N
\varepsilon_1\int_B|\nabla^2\phi|^2\xi^2+\frac{C_N}{\varepsilon_1}\int_B|d\phi|^2|\nabla\xi|^2,
\end{equation}
\begin{equation}
\label{d20} II\leq C_N
\varepsilon_1\int_B|\nabla^2\phi|^2\xi^2+\frac{C_N}{\varepsilon_1}\int_B|d\phi|^4\xi^2.
\end{equation}
Choosing $\eta=|d\phi|\xi$ in (\ref{d13}), we obtain
\begin{equation}
\label{d21} \||\nabla\psi||d\phi|\xi\|_{L^{\frac 4
3}(B)}+\|\psi|d\phi|\xi\|_{L^4(B)}\leq
C\|\nabla(|d\phi|\xi)|\psi|\|_{L^{\frac 4 3}(B)}.
\end{equation}
Because
\begin{eqnarray}
\label{d22} \|\nabla(|d\phi|\xi)|\psi|\|^2_{L^{\frac 4
3}(B)}&\leq&2\||\nabla^2\phi|\psi\xi\|^2_{L^{\frac 4
3}(B)}+2\||d\phi||\nabla\xi||\psi|\|^2_{L^{\frac 4 3}(B)}\nonumber
\\
&\leq&2(\int_B|\nabla^2\phi|^2\xi^2(\int_B|\psi|^4)^{\frac 1
2}+2(\int_B|d\phi|^2|\nabla\xi|^2)(\int_B|\psi|^4)^{\frac 1
2}\nonumber \\
&=&2(\int_B|\psi|^4)^{\frac 1
2}(\int_B|\nabla^2\phi|^2\xi^2+\int_B|d\phi|^2|\nabla\xi|^2),
\end{eqnarray}
choose $R>0$ small enough such that
\begin{equation}
\label{d23} C^2\max \{2C_N,C_N^2\}\int_{D(x_0,R)}|\psi|^4<\frac 1
8 ,
%\qquad C^2C_N^2\int_{D(x_0,R)}|\psi|^4<\frac 1 8 ,
\end{equation}
we get
\begin{eqnarray}
\label{d24} IV&\leq&C_N\||\nabla\psi||d\phi|\xi\|_{L^{\frac 4
3}(B)}\|\psi|d\phi|\xi\|_{L^4(B)}\nonumber \\
&\leq&C_NC^2\|\nabla(|d\phi|\xi)|\psi|\|^2_{L^{\frac 4 3}(B)}
\quad (\hbox{ by }(\ref{d21}))\nonumber \\ &\leq&
2C_NC^2(\int_B|\psi|^4)^{\frac 1
2}(\int_B|\nabla^2\phi|^2\xi^2+\int_B|d\phi|^2|\nabla\xi|^2)\quad
({\rm by (\ref{d22})})\nonumber
\\
&<& \frac 1 8 \int_B|\nabla^2\phi|^2\xi^2+ \frac 1 8
\int_B|d\phi|^2|\nabla\xi|^2.
\end{eqnarray}
Similarly, we can estimate the terms $V$ and $VI$ as follows.
\begin{eqnarray}
\label{d25} V&\leq& C_N
\varepsilon_1\int_B|d\phi|^2\xi^2|\psi|^4+\frac{C_N}{\varepsilon_1}\int_B|d\phi|^4\xi^2\nonumber
\\
&\leq& C_N \varepsilon_1(\int_B|\psi|^4)^{\frac 1
2}(\int_B|\psi|^4|d\phi|^4\xi^4)^{\frac 1 2}+\frac{C_N}{\varepsilon_1}\int_B|d\phi|^4\xi^2\nonumber \\
&\leq&C_N \varepsilon_1(\int_B|\psi|^4)^{\frac 1
2}C^2|\nabla(|d\phi|\xi)|\psi||^2_{L^{\frac 4
3}(B)}+\frac{C_N}{\varepsilon_1}
\int_B|d\phi|^4\xi^2\quad ({\rm by (\ref{d21})})\nonumber \\
&\leq&
2C_NC^2\varepsilon_1(\int_B|\psi|^4)(\int_B|\nabla^2\phi|^2\xi^2+\int_B|d\phi|^2|\nabla\xi|^2)+\frac{C_N}{\varepsilon_1}
\int_B|d\phi|^4\xi^2\nonumber \\
&<&\frac 1 8 \int_B|\nabla^2\phi|^2\xi^2+ \frac 1 8
\int_B|d\phi|^2|\nabla\xi|^2+\frac{C_N}{\varepsilon_1}
\int_B|d\phi|^4\xi^2,
\end{eqnarray}
and
\begin{eqnarray}
\label{d26} VI&=&
C_N\int_B(|\nabla^2\phi|\xi)(|\psi|^2|d\phi|\xi)\nonumber\\
&\leq&\frac 1 2 \int_B|\nabla^2\phi|^2\xi^2+\frac 1 2
C_N^2\int_B|\psi|^4|d\phi|^2\xi^2\nonumber\\
&\leq&\frac 1 2 \int_B|\nabla^2\phi|^2\xi^2+\frac 1 2
C_N^2(\int_B|\psi|^4)^{\frac 1
2}(\int_B|\psi|^4|d\phi|^4\xi^4)^{\frac 1 2}\nonumber\\
&\leq&\frac 1 2 \int_B|\nabla^2\phi|^2\xi^2+\frac 1 2
C_N^2(\int_B|\psi|^4)^{\frac 1
2}C^2|\nabla(|d\phi|\xi)|\psi||^2_{L^{\frac 4 3}(B)}  \nonumber\\
&\leq&\frac 1 2
\int_B|\nabla^2\phi|^2\xi^2+C_N^2C^2(\int_B|\psi|^4)(\int_B|\nabla^2\phi|^2\xi^2+
\int_B|d\phi|^2|\nabla\xi|^2)  \nonumber\\
&<&\frac 1 2 \int_B|\nabla^2\phi|^2\xi^2+\frac 1 8
\int_B|\nabla^2\phi|^2\xi^2+ \frac 1 8
\int_B|d\phi|^2|\nabla\xi|^2.
\end{eqnarray}
Putting (\ref{d19}), (\ref{d20}), (\ref{d24}), (\ref{d25}) and
(\ref{d26}) into (\ref{d18}) gives
\begin{equation}
\label{d27} \int_{D(x_0,R)}|\nabla^2\phi|^2\xi^2\leq
C(\int_{D(x_0,R)}|d\phi|^2|\nabla\xi|^2+\int_{D(x_0,R)}|d\phi|^4\xi^2).
\end{equation}
Now for $\varepsilon>0$, let $\rho>0$ be as in Lemma 2.4, and with
$D(x_1,\rho)\subset D(x_0,R)$,  choose a cut-off function $\xi\in
C_0^\infty(D(x_1,\rho)), $  $0\leq \xi \leq 1$ such that
$$\xi\equiv 1 \quad {\rm in } \quad D(x_1,\frac{\rho}{2});\quad
|\nabla\xi|\leq \frac 4 \rho \quad {\rm in } \quad D(x_1,\rho),$$
denote $B_\rho := D(x_1,\rho)$ for simplicity, one derives
\begin{eqnarray*}
\int_{B_\rho}|d\phi|^4\xi^2&=&\int_{B_\rho}|d\phi|^2(|d\phi|\xi)^2\nonumber
\\
&\leq&\varepsilon \int_{B_\rho}|\nabla(|d\phi|\xi)|^2+C\varepsilon
(\int_{B_\rho}|\psi|^4|d\phi|^4\xi^4)^{\frac 1 2} \quad ({\rm
by\quad
  Lemma 2.4}) \nonumber \\
&\leq&\varepsilon \int_{B_\rho}|\nabla^2\phi|^2\xi^2+\varepsilon
\int_{B_\rho}|d\phi|^2|\nabla\xi|^2+C\varepsilon
(\int_{B_\rho}|\psi|^4|d\phi|^4\xi^4)^{\frac 1 2},
\end{eqnarray*}
it follows from (\ref{d21}) that
$$\||\psi||d\phi|\xi\|_{L^4(B_\rho)}\leq
C\|\nabla(|d\phi|\xi)|\psi|\|_{L^{\frac 4 3}(B_\rho)},$$ then, by
an argument similar to the one used in the proof of (\ref{d22}),
we can get $$\||\psi||d\phi|\xi\|^2_{L^4(B_\rho)}\leq 2C^2
(\int_{B_\rho}|\psi|^4)^{\frac 1
2}(\int_{B_\rho}|\nabla^2\phi|^2\xi^2+\int_{B_\rho}|d\phi|^2|\nabla\xi|^2),$$
thus,
\begin{eqnarray}
\label{d28} \int_{B_\rho}|d\phi|^4\xi^2&\leq&\varepsilon
\int_{B_\rho}|\nabla^2\phi|^2\xi^2+\varepsilon
\int_{B_\rho}|d\phi|^2|\nabla\xi|^2\nonumber \\
&&+C'\varepsilon
(\int_{B_\rho}|\nabla^2\phi|^2\xi^2+\int_{B_\rho}|d\phi|^2|\nabla\xi|^2)\nonumber
\\
&=& C''\varepsilon (\int_{B_\rho}|\nabla^2\phi|^2\xi^2+
\int_{B_\rho}|d\phi|^2|\nabla\xi|^2).
\end{eqnarray}
On the other hand, by (\ref{d27}), we have
\begin{equation}
\label{d29} \int_{B_\rho}|\nabla^2\phi|^2\xi^2\leq
C(\int_{B_\rho}|d\phi|^2|\nabla\xi|^2+\int_{B_\rho}|d\phi|^4\xi^2),
\end{equation}
substituting (\ref{d28}) into (\ref{d29}) yields
$$\int_{B_\rho}|\nabla^2\phi|^2\xi^2\leq
C\int_{B_\rho}|d\phi|^2|\nabla\xi|^2,$$ hence,
\begin{equation}
\label{d30} \int_{D(x_1,\frac{\rho}{2})}|\nabla^2\phi|^2\leq
\frac{C}{\rho^2}\int_{D(x_1,\rho)}|d\phi|^2.
\end{equation}
Covering $D(x_0,\frac{R}{2})$ with $\{D(x_1,\frac{\rho}{2})\}$ and
using (\ref{d30}) we obtain (\ref{d8}). \hfill $\Box$ \par Now  we
are in the position to give the \vskip12pt \noindent {\it Proof of
Theorem 2.3.} First, we show that $\phi\in W^{2,2}\cap
W^{1,4}(D(x_0,\frac R 2 ),N)$. This can be done just by replacing
weak derivatives by difference quotients in the proof of Lemma
2.5. Denote
$$\Delta^h_i\phi(x):=\frac{\phi(x+hE_i)-\phi(x)}{h}, \quad ({\rm
if \quad dist }(x,\partial D)>|h|),$$ where $(E_1,E_2,\cdots,
E_K)$ is an orthonarmal basis of $\RR^K$, $h\in \RR$.
$\Delta^h:=(\Delta^h_1,\Delta^h_2,\cdots,\Delta^h_K)$. Let
$\zeta:=\Delta^{-h}_\gamma (\xi^2\Delta^h_\gamma\phi)$, then,
similar to (\ref{d27}), we have
\begin{equation}
\label{d31} \int_{D(x_0,R)}|\nabla(\Delta^h\phi)|^2\xi^2\leq C
\int_{D(x_0,R)}|\Delta^h\phi|^2|\nabla\xi|^2+C\int_{D(x_0,R)}
|\nabla\phi|^2|\Delta^h\phi|^2\xi^2.
\end{equation}
Since (cf. \cite{J1}, p.382) $$C
\int_{D(x_0,R)}|\Delta^h\phi|^2|\nabla\xi|^2\leq C
\int_{D(x_0,R)}|\nabla\phi|^2|\nabla\xi|^2,
$$
applying Lemma 2.4 to the right hand side, we obtain the following
estimate analogous to (\ref{d30}):
$$\int_{D(x_1,\frac \rho 2 )}|\nabla(\Delta^h\phi)|^2\xi^2\leq
\frac{C}{\rho^2}\int_{D(x_1,\rho)}|d\phi|^2,$$ from which it
follows that the weak derivative $\nabla^2\phi$ exists and
(\ref{d8}) still holds true with $\rho$ sufficiently small and $C_1>0$ which depends on
$|\phi|_{C^0(D,N)}$ and $\rho$.
\par
Next, since $\phi\in W^{2,2}$, we have that $\phi\in W^{1,p}$ for
any $p>0$, thus, the right hand side of (\ref{d2}) is in $L^p
(p>2)$, so $\psi \in C^{0,\gamma}$ for some $\gamma >0$. By the
elliptic estimates for the equation (\ref{d1}), we have $\phi \in W^{2,p}$ for any $p>2$,
thus $\phi \in C^{1,\gamma}$. By the elliptic estimates for the equation (\ref{d2}), we have
$\psi \in C^{1,\gamma}$. Then the standard arguments yield that both $\phi$ and $\psi$ are smooth.
 This completes the proof.
 \hfill $\Box$
\section{ Energy identities}
\addtocounter{equation}{-51} First, using the elliptic estimates,
we can establish the following vanishing theorem which will be
used later in obtaining the energy identities, and in which we see
that the $W^{1,2}$-norm of $\phi$ and $L^4$-norm of $\psi$ play an
important role in the analytic properties of $\phi$ and $\psi$. We
note that these two norms are conformally invariant.
\par
\vskip12pt
 \noindent
 {\bf Theorem 3.1.} {\it Let $(M^2,h_{\alpha\beta})$ be a compact Riemann surface
 with a fixed spin structure, and $(N^n,g_{ij})$ be a compact Riemannian manifold.
 There is a constant $\varepsilon_0 >0$ small
 enough such that if $(\phi,\psi)$ is a smooth solution of (\ref{b1}) and
 (\ref{b2}) satisfying
 \begin{equation}
 \label{e4.1}
 \int_M(|d\phi|^2+|\psi|^4) <\varepsilon_0,
 \end{equation}
then $\phi$ is constant and consequently $\partial \hskip-2.2mm
\slash \psi^i\equiv 0,$ $i=1,2,\cdots,n.$ }
\par
\vskip12pt
 \noindent {\it Proof.}
In the sequel, we write $\|\cdot\|_{D,k,p}$ for the $L^{k,p}$-norm
on the domain $D$, and if there is no confusion, we may drop the
subscript $D$.
 Embed $N$ into some $\R^K$ isometrically,
then from the $\phi$-equation (\ref{b9}) we have
$$|\Delta\phi|\leq \|A\|_\infty|d\phi|^2+C_N|d\phi||\psi|^2,$$
where $C_N>0$ is a constant depending only on $N$,
$\|A\|_\infty:=\max_N|A|.$ It follows from the above inequality
that
\begin{eqnarray*} %\begin{array}\label{e4.2}
\|\Delta\phi\|_{0,\frac{4}{3}} & \leq&
\|A\|_\infty\||d\phi|^2\|_{0,\frac{4}{3}}+C_N\||d\phi||\psi|^2\|_{0,\frac{4}{3}}\\
&\le& C(\|d\phi\|^2_{0,2}+\|\psi\|^2_{0,4})\|d\phi\|_{0,4}\\ &\le&
C(\|d\phi\|^2_{0,2}+\|\psi\|^2_{0,4})\|d\phi\|_{1,4/3}.
\end{eqnarray*}
%Because
% \begin{eqnarray*}
%(|d\phi|^2)_{0,\frac{4}{3}}&\leq&|d\phi|_{0,2}|d\phi|_{0,4}\\
%&\leq&C|\phi|_{2,\frac{4}{3}}|d\phi|_{0,2},
%\end{eqnarray*}
%and
%\begin{eqnarray*}
%(|d\phi||\psi|^2)_{0,\frac{4}{3}}&\leq&|\psi|^2_{0,4}|d\phi|_{0,4}\\
%&\leq& C|\psi|^2_{0,4}|\phi|_{2,\frac{4}{3}},
%\end{eqnarray*}
%we have
%\begin{eqnarray*}
%|\Delta\phi|_{0,\frac{4}{3}} \leq |A|_\infty
%C|\phi|_{2,\frac{4}{3}}|d\phi|_{0,2}+C_NC|\psi|^2_{0,4}|\phi|_{2,\frac{4}{3}}.
%\end{eqnarray*}
%Hence
%$$|\phi|_{2,\frac{4}{3}}\leq C(|A|_\infty|d\phi|_{0,2}+|\psi|^2_{0,4})|\phi|_{2,\frac{4}{3}}.$$
%If
%$
If $\int_M(|d\phi|^2+|\psi|^4) <\varepsilon_0$ for small
$\varepsilon_0>0$, then $\phi\equiv const$. From equation
(\ref{b8}), we have $\partial \hskip-2.2mm \slash \psi^i\equiv 0,$
$i=1,2,\cdots,n.$ \hfill $\Box$
\par
Now we prove the small energy regularity. Since the problem is
local, we assume that $M$ is flat. \vskip12pt \noindent {\bf
Theorem 3.2 ($\varepsilon-$regularity theorem).} {\it There is an
$\varepsilon_0>0$ such that if $(\phi,\psi):
(D,\delta_{\alpha\beta})\to (N,g_{ij})$ is a $C^\infty$
Dirac-harmonic map satisfying
\begin{equation}
\label{c1} \int_D(|d\phi|^2+|\psi|^4)<\varepsilon_0,
\end{equation}
then
\begin{equation}
\label{5.8} \|d\phi\|_{\widetilde{D},1,p}\leq
C(\widetilde{D},p)\|d\phi\|_{D,0,2},
\end{equation}
\begin{equation}
\label{5.11} \|\nabla\psi\|_{\widetilde{D},1,p}\leq
C(\widetilde{D},p)\|\psi\|_{D,0,4},
\end{equation}
\begin{equation}
\label{5.12} \|\nabla\psi\|_{C^0(\widetilde{D})}\leq
C(\widetilde{D})\|\psi\|_{D,0,4},\quad
\|\psi\|_{C^0(\widetilde{D})} \leq
C(\widetilde{D})\|\psi\|_{D,0,4}
\end{equation}
 $\forall \widetilde{D}\subset D, p>1,$
where $C(\widetilde{D},p)>1$ is a constant depending only on
$\widetilde{D}$ and $p$. } \vskip12pt To prove this theorem, we
first estimate $|d\phi|$. \vskip12pt \noindent
 {\bf Lemma 3.3} {\it There is an
$\varepsilon_0>0$ such that if $(\phi,\psi):
(D,\delta_{\alpha\beta})\to (N,g_{ij})$ is a $C^\infty$
Dirac-harmonic map satisfying (\ref{c1}), then
\begin{equation}
\label{5.3} \|\phi\|_{D^1,1,4}\leq
C(D^1)\sqrt{\varepsilon_0},\quad \forall D^1 \quad {\rm with}\quad
\overline{D^1}\subset D,
\end{equation}
where $C(D^1)>0$ is a constant depending only on $D^1$. }
\vskip12pt \noindent {\it Proof.}
 Choose a cut-off function $\eta: 0\leq \eta \leq 1, $ with $\eta|_{D^1}\equiv 1$
 and ${\rm Supp}\eta\subset D$. By (\ref{b9}) we have
\begin{eqnarray*}
|\Delta(\eta\phi)|&\leq&
C(|\phi|+|d\phi|)+\|A\|_\infty|d\phi|(|d(\eta\phi)|+|\phi d\eta|)
+|\eta\alpha| \\ &\leq&
\|A\|_\infty|d\phi||d(\eta\phi)|+C(|\phi|+|d\phi|)+|\eta\alpha|,
\end{eqnarray*}
where $\alpha :=P(\cal
A(d\phi(e_\alpha),e_\alpha\cdot\psi);\psi),$ thus, for any $p>1$,
\begin{equation}
\label{5.4} \|\Delta(\eta\phi)\|_{0,p}\leq
\|A\|_\infty\||d\phi||d(\eta\phi)|\|_{0,p}+C\|\phi\|_{1,p}+\|\eta\alpha\|_{0,p}.
\end{equation}
Let $p=\frac{4}{3}$, and without loss of generality we assume
$\int_D\phi=0$ so that $\|\phi\|_{1,p}\leq C'\|d\phi\|_{0,p},$
then $$\|A\|_\infty\||d\phi||d(\eta\phi)|\|_{0,\frac{4}{3}}\leq
\|A\|_\infty \|\eta\phi\|_{1,4}\|d\phi\|_{0,2},$$ from this and
(\ref{5.4}) we have $$ \|\eta\phi\|_{2,\frac{4}{3}}\leq
C(\|A\|_\infty
\|\eta\phi\|_{1,4}\|d\phi\|_{0,2}+\|d\phi\|_{0,\frac{4}{3}}+\|\eta\alpha\|_{0,\frac{4}{3}}).
$$ Bythe Sobolev inequality, $\|\eta\phi\|_{1,4}\leq
C'\|\eta\phi\|_{2,\frac{4}{3}}$, so,
\begin{equation}
\label{5.5} (C^{'-1}-C\|A\|_\infty
\|d\phi\|_{0,2})\|\eta\phi\|_{1,4}\leq
C(\|d\phi\|_{0,\frac{4}{3}}+\|\eta\alpha\|_{0,\frac{4}{3}}),
\end{equation}
moreover,
\begin{eqnarray*}
\|\eta\alpha\|_{0,\frac{4}{3}}&\leq& C_N\||\psi|^2|\eta
d\phi|\|_{0,\frac{4}{3}}\\
       &=& C_N\||\psi|^2|d(\eta \phi)-\phi d\eta|\|_{0,\frac{4}{3}}\\
&\leq& C\||\psi|^2|d(\eta
\phi)|\|_{0,\frac{4}{3}}+C\||\psi|^2|\phi||
d\eta|\|_{0,\frac{4}{3}}\\ &\leq&
C\|\psi\|^2_{0,4}\|\eta\phi\|_{1,4}+C\|\psi\|^2_{0,\frac{4}{3}}\\
&\leq& C\|\psi\|^2_{0,4}\|\eta\phi\|_{1,4}+C\|\psi\|^2_{0,4},
\end{eqnarray*}
putting this into (\ref{5.5}) we get: $$\|\eta\phi\|_{1,4}\leq
C(\|d\phi\|_{0,\frac{4}{3}}+\sqrt{\varepsilon_0}\|\eta\phi\|_{1,4}
+\|\psi\|^2_{0,4}),$$ which yields $$\|\eta\phi\|_{1,4}\leq
C(\|d\phi\|_{0,\frac{4}{3}}+\|\psi\|^2_{0,4})<2\sqrt{\varepsilon_0}C.$$
\hfill $\Box$ \par Next we estimate $\psi$.
 \vskip12pt \noindent {\bf
Lemma 3.4.} {\it There is an $\varepsilon_0>0$ such that if
$(\phi,\psi): (D,\delta_{\alpha\beta})\to (N,g_{ij})$ is a
$C^\infty$ Dirac-harmonic map satisfying (\ref{c1}), then
\begin{equation}
\label{5.6} \|\psi\|_{D^2,0,q}\leq C(D^2)\|\psi\|_{D,0,4},\quad
\forall q>1, \quad \forall D^2\quad {\rm with}\quad
\overline{D^2}\subset D,
\end{equation}
where $C(D^2)>0$ is a constant depending only on $D^2$. }
\vskip12pt \noindent {\it Proof.} Choose a cut-off function $\eta:
0\leq \eta \leq 1, $ with $\eta|_{D^2}\equiv 1$ and ${\rm
Supp}\eta\subset D$. For any $i=1,2,\cdots, n$, $\psi^i$ is an
ordinary spinor field, and $\xi^i:=\eta\psi^i$ has compact support
in $D$, so, by the well-known Lichnerowitz's formula, we have $$
\partial \hskip -2.2mm \slash^2 \xi^i=-\Delta
\xi^i+\frac{1}{4}R\xi^i=-\Delta \xi^i$$ because the scalar
curvature $R\equiv 0$ on $D$. Integrating this yields
\begin{eqnarray*}
\int_{D}|\nabla\xi^i|^2 &=& \int_{D}|\partial \hskip -2.2mm \slash
\xi^i|^2 \\ &=& \int_{D}|\partial \hskip -2.2mm \slash (\eta
\psi^i)|^2 \\ &=& \int_{D}|\nabla\eta\cdot \psi^i+\eta\partial
\hskip -2.2mm \slash \psi^i|^2 \\ &\leq & C(\int_{D}|\psi^i|^2+
\int_{D}|\partial \hskip -2.2mm \slash \psi^i|^2)\\ &\leq &
C(\int_{D}|\psi^i|^2+ C\int_{D}|d\phi|^2|\psi|^2),
\end{eqnarray*}
hence,
\begin{eqnarray*}
\|\nabla\xi^i\|_{D,0,2} &\leq&
C(\|\psi^i\|_{D,0,2}+\|d\phi\|^2_{D,0,4}\|\psi\|^2_{D,0,4})\\
&\leq& C\|\psi\|_{D,0,4}(1+\|d\phi\|_{D,0,4})\\ &\leq&
C'\|\psi\|_{D,0,4},
\end{eqnarray*}
from which it follows that
$$\|\psi\|_{D^2,0,q}<C\|\psi\|_{D,0,4}.$$ \hfill $\Box$ \vskip12pt
\noindent {\bf Lemma 3.5.} {\it There is an $\varepsilon_0>0$ such
that if $(\phi,\psi): (D,\delta_{\alpha\beta})\to (N,g_{ij})$ is a
$C^\infty$ Dirac-harmonic map satisfying (\ref{c1}), then
\begin{equation}
\label{5.7} \|d\phi\|_{D^2,0,4}\leq C(D^2)\|d\phi\|_{D,0,2}, \quad
\forall D^2\quad {\rm with}\quad \overline{D^2}\subset D,
\end{equation}
where $C(D^2)>0$ is a constant depending only on $D^2$. }
\vskip12pt \noindent {\it Proof.}
 Choose a cut-off function $\eta: 0\leq \eta \leq 1, $ with $\eta|_{D^2}\equiv 1$ and ${\rm Supp}\eta\subset D$.
By (\ref{5.5}), we have
\begin{eqnarray*}
\|\eta\phi\|_{D,1,4}&\leq &
C(\|d\phi\|_{D,0,\frac{4}{3}}+\|\eta\alpha\|_{D,0,\frac{4}{3}})\\
&\leq & C(\|d\phi\|_{D,0,2}+\|\eta\alpha\|_{D,0,\frac{4}{3}}),
\end{eqnarray*}
and it is clear that
\begin{eqnarray*}
\|\alpha\|_{D,0,\frac{4}{3}}&\leq &
C\||\psi|^2|d\phi|\|_{D,0,\frac{4}{3}}\\ &\leq &
C(\int_{D}|\psi|^8)^{\frac{1}{4}}(\int_{D}|d\phi|^2)^{\frac{1}{2}}\\
&\leq & C\|\psi\|_{D,0,4}\|d\phi\|_{D,0,2},
\end{eqnarray*}
therefore, $$\|\eta\phi\|_{D,1,4}\leq  C\|d\phi\|_{D,0,2}.$$
\hfill $\Box$

\par
\vskip12pt \noindent {\it Proof of Theorem 3.2.}  Choose
$\widetilde{D}\subset D^3\subset D^2\subset D^1\subset D$. Also
choose a cut-off function
 $\eta: 0\leq \eta \leq 1, $ with $\eta|_{D^3}\equiv 1$ and ${\rm Supp}\eta\subset D^2$.
By (\ref{5.4}) on $D^2$ for $p=2$ (we temporarily assume
$\int_{D^2}\phi =0$):
\begin{eqnarray*}
\|\eta\phi\|_{D^2,2,2} &\leq &
C[\|A\|_\infty\|d(\eta\phi)\|_{D^2,0,4}\|d\phi\|_{D^2,0,4}+\|\phi\|_{D^2,1,2}
       +\||\psi|^2|d\phi|\|_{D^2,0,2}]\\
&\leq &
C(\|A\|_\infty\|\eta\phi\|_{D^2,1,4}\|d\phi\|_{D^2,0,4}+\|\phi\|_{D^2,1,2}
       +\|\psi\|^2_{D^2,0,8}\|d\phi\|_{D^2,0,4}).
\end{eqnarray*}
Since $$ \|\eta\phi\|_{D^2,1,4}\leq
C(D^2)\|\eta\phi\|_{D^2,2,\frac{4}{3}}\leq
C\|\eta\phi\|_{D^2,2,2}, $$ using this and (\ref{5.3}) we otain $$
\|A\|_\infty\|\eta\phi\|_{D^2,1,4}\|d\phi\|_{D^2,0,4}\leq
C'\sqrt{\varepsilon_0}\|A\|_\infty\|\eta\phi\|_{D^2,2,2},$$ which
yields
\begin{eqnarray*}
(1-CC'\sqrt{\varepsilon_0}\|A\|_\infty)\|\eta\phi\|_{D^2,2,2}
&\leq & C(\|\phi\|_{D^2,1,2}+\|\psi\|^2_{D^2,0,8}
\|d\phi\|_{D^2,0,4})
\\
 &\leq & C\|\phi\|_{D^2,1,4},
\end{eqnarray*}
therefore, $$\|\phi\|_{D^2,2,2}\leq C\|\phi\|_{D^2,1,4} \leq
\|d\phi\|_{D^2,0,4}.$$ By the Sobolev inequality, we have
\begin{equation}
\label{5.9} \|d\phi\|_{D^3,0,p}<C\|d\phi\|_{D^2,0,4}, \quad
\forall p>1.
\end{equation}
This also holds for $\phi$ without $\int_{D^2}\phi=0.$
\par
Now for $\widetilde{D}\subset D^3$, as above, we choose a cut-off
function $\eta: 0\leq \eta \leq 1, $ with
$\eta|_{\widetilde{D}}\equiv 1$ and ${\rm Supp}\eta\subset D^3$.
By (\ref{5.4}) on $D^3$ for any $p>1$ (we again temporarily assume
$\int_{D^3}\phi =0$) we have: $$\|\eta\phi\|_{D^3,2,p}\leq C[
\|A\|_\infty
\||d\phi||d(\eta\phi)|\|_{D^3,0,p}+\|\phi\|_{D^3,1,p}+\||\psi|^2|d\phi|\|_{D^3,0,p}].
$$ Using (\ref{5.9}), we obtain:
$$\||d\phi||d(\eta\phi)|\|_{D^3,0,p}\leq
(\int_{D^3}|d\phi|^{2p})^{\frac{1}{p}}\leq C\|d\phi\|^2_{D^2,0,4}
\leq C\|d\phi\|_{D^2,0,4}$$ and
\begin{eqnarray*}
\||\psi|^2|d\phi|\|_{D^3,0,p}&=&(\int_{D^3}|\psi|^{2p}|d\phi|^p)^{\frac{1}{p}}\\
&\leq & (\|\psi\|_{D^3,0,4p})^2\|d\phi\|_{D^3,0,2p}\\ &\leq &
C\|d\phi\|_{D^2,0,4},
\end{eqnarray*}
we conclude that $$ \|\eta\phi\|_{D^3,2,p}\leq
C\|d\phi\|_{D^2,0,4},$$ from which and (\ref{5.7}) it follows that
$$ \|\phi\|_{\widetilde{D},2,p}\leq C\|d\phi\|_{D^2,0,4} \\ \leq
C\|d\phi\|_{D,0,2}. $$ Clearly, this also holds for $\phi$ without
$\int_{D^3}\phi =0.$ Now let us give the estimates for $\psi$. By
Lemmas 3.4 and 3.5, for $D^1\subset D$, the following estimates
hold
\begin{equation}
\label{5.13} \|\psi\|_{D^1,0,p}\leq C\|\psi\|_{D,0,4},
\end{equation}
\begin{equation}
\label{5.14} \|d\phi\|_{D^1,1,p}\leq C\|d\phi\|_{D,0,2}.
\end{equation}
On the other hand, for $D^1\subset D$, by an argument similar to
the one used in the proof of Lemma 3.4, we have
\begin{equation}
\label{5.15} \|\nabla\psi^i\|_{D^1,0,2}\leq C\|\psi\|_{D,0,4}.
\end{equation}
Computing directly, one gets
$$\D^2\psi^i=-(\wa\wa\psi)^i+\frac{1}{2}R^i_{kpj}\nabla\phi^p\cdot\nabla\phi^j\cdot\psi^k,$$
from which we have
\begin{eqnarray}
\label{5.16}
&&\Delta\psi^i+(\Gamma^i_{jk,p}+\Gamma^i_{pm}\Gamma^m_{jk})\phi^p_\alpha\phi^j_\alpha\psi^k -\frac{1}{2}R^i_{kpj}\nabla\phi^p\cdot\nabla\phi^j\cdot\psi^k  \nonumber \\
&&
\hskip108pt+\Gamma^i_{jk}\phi^j_{\alpha\alpha}\psi^k+2\Gamma^i_{jk}\phi^j_\alpha\nabla_{e_\alpha}\psi^k=0.
\end{eqnarray}
For any $\eta\in C^{\infty}(D,R)$ with $0\leq\eta\leq 1$, we have
\begin{equation}
\label{5.17} |\Delta (\eta\psi)|\leq
C(|\psi|+|\nabla\psi|+|d\phi|^2|\psi|+|\nabla^2\phi||\psi|+|d\phi||\nabla\psi|),
\end{equation}
on $D^2\subset D^1$. Choose a cut-off function $\eta: 0\leq \eta
\leq 1, $ with $\eta|_{D^2}\equiv 1$ and ${\rm Supp}\eta\subset
D^1$. For any $p>1$ we have
\begin{eqnarray}
\label{5.18} \|\eta\psi\|_{D^1,2,p} &\leq &
C_p[\|\psi\|_{D^1,0,p}+\|\nabla\psi\|_{D^1,0,p}+\||d\phi|^2|\psi|\|_{D^1,0,p}\nonumber
\\
&&+\||\nabla^2\phi||\psi|\|_{D^1,0,p}+\||d\phi||\nabla\psi|\|_{D^1,0,p}].
\end{eqnarray}
Let $p=\frac{4}{3}$,  by (\ref{5.13}), (\ref{5.14}) and
(\ref{5.15}), we obtain
\begin{eqnarray*}
\|\psi\|_{D^2,2,\frac{4}{3}} &\leq &
C[\|\psi\|_{D^1,0,\frac{4}{3}}+\|\nabla\psi\|_{D^1,0,\frac{4}{3}}+\|d\phi\|^2_{D^1,0,4}\|\psi\|_{D^1,0,4}\nonumber
\\
&&+\|\nabla^2\phi\|_{D^1,0,2}\|\psi\|_{D^1,0,4}+\|d\phi\|_{D^1,0,4}\|\nabla\psi\|_{D^1,0,2}]\nonumber
\\ &\leq & C\|\psi\|_{D,0,4},
\end{eqnarray*}
from which it follows that
\begin{equation}
\label{5.19} \|\psi\|_{D^2,1,4}\leq C\|\psi\|_{D,0,4},
\end{equation}
and consequently, $$\|\psi\|_{C^0(D^2)}\leq C\|\psi\|_{D,0,4},$$
This proves the second inequality in (\ref{5.12}).
\par
Now for $D^3\subset D^2$, choose a cut-off function $\eta: 0\leq
\eta \leq 1, $ with $\eta|_{D^3}\equiv 1$ and ${\rm
Supp}\eta\subset D^2$. By (\ref{5.18}) on $D^2$ for $p=2$ we have:
\begin{eqnarray*}
\|\eta\psi\|_{D^2,2,2} &\leq &
C[\|\psi\|_{D^2,0,2}+\|\nabla\psi\|_{D^2,0,2}+\||d\phi|^2|\psi|\|_{D^2,0,2}\nonumber
\\
&&+\||\nabla^2\phi||\psi|\|_{D^2,0,2}+\||d\phi||\nabla\psi|\|_{D^2,0,2}]\nonumber
\\ &\leq &
C[\|\psi\|_{D^2,0,4}+\|\psi\|_{D,0,4}+\|d\phi\|^2_{D^2,0,8}\|\psi\|_{D^2,0,4}\nonumber
\\
&&+\|\nabla^2\phi\|_{D^2,0,4}\|\psi\|_{D^2,0,4}+\|d\phi\|_{D^2,0,4}\|\nabla\psi\|_{D^2,0,4}]\nonumber
\\ &\leq & C\|\psi\|_{D,0,4},
\end{eqnarray*}
consequently,
\begin{equation}
\label{5.20} \|\psi\|_{D^3,1,p}\leq C\|\psi\|_{D,0,4},\quad
\forall p>1.
\end{equation}
We again use (\ref{5.18}) on $D^3$ for a cut-off function $\eta:
0\leq \eta \leq 1, $ with $\eta|_{\widetilde{D}} \equiv 1$ and
${\rm Supp}\eta\subset D^3$. By (\ref{5.13}) (\ref{5.14})
(\ref{5.15}) and (\ref{5.20}) we have:
$$\|\eta\psi\|_{D^3,2,p}\leq C\|\psi\|_{D,0,4},\quad \forall p>1.
$$ that is, $$\|\nabla\psi\|_{\widetilde{D},1,p}\leq
C(\widetilde{D},p)\|\psi\|_{D,0,4}.$$ We therefore obtain
$$\|\nabla\psi\|_{C^0(\widetilde{D})} \leq
C(\widetilde{D})\|\psi\|_{D,0,4}. $$ This proves the first
inequality in (\ref{5.12}). \hfill $\Box$ \vskip12pt
 \noindent
{\bf Theorem 3.6 (Energy identities).}  {\it Let
$\{(\phi_k,\psi_k): M\to \S^n\}$ be a sequence of smooth
Dirac-harmonic maps with uniform bounded energy:
$$E(\phi_k,\psi_k)\leq\Lambda<+\infty,$$ and assume that
$\{(\phi_k,\psi_k)\}$ weakly converges to a Dirac-harmonic map
$(\phi,\psi)$ in $W^{1,2}(M,\S^n)\times L^4(\Sigma M\otimes
\RR^{n+1})$, then the blow-up set $$S:=\cap_{r>0}\{x\in M|{\rm lim
inf}_{k\to
\infty}\int_{Dx,r)}(|d\phi_k|^2+|\psi_k|^4)\geq\varepsilon_0\}$$
is a finite set of points $\{p_1,p_2,\cdots,p_I\}$, and there are
a subsequence, still denoted by $\{(\phi_k,\psi_k)\}$, and a
finite set of Dirac-harmonic maps $(\sigma_i^l,\xi_i^l): \S^2\to
\S^n, i=1,2,\cdots,I; $ $l=1,2,\cdots,L_i$ such that
\begin{equation}
\label{c4} {\rm lim }_{k\to
\infty}E(\phi_k)=E(\phi)+\Sigma^I_{i=1}\Sigma_{l=1}^{L_i}E(\sigma^l_i),
\end{equation}
\begin{equation}
\label{c5} {\rm lim }_{k\to
\infty}E(\psi_k)=E(\psi)+\Sigma^I_{i=1}\Sigma_{l=1}^{L_i}E(\xi^l_i),
\end{equation}
where $\varepsilon_0$ is as in Theorem 3.2,
$E(\phi_k):=\int_M|d\phi_k|^2,$ $E(\psi_k):=\int_M|\psi_k|^4,
E(\phi_k,\psi_k):=\int_M(|d\phi_k|^2+|\psi_k|^4).$
 }
\vskip12pt \noindent {\it Proof.} Since
$E(\phi_k,\psi_k)\leq\Lambda<+\infty,$ the blow-up set $S$ must be
finite. So we can find small disks $D_{\delta_i}$ for each $p_i$
such that $D_{\delta_i}\cap D_{\delta_j}=\emptyset$ for $i\not= j,$ $
i,j=1,2,\cdots,I.$  On $M\setminus\cup_{i=1}^ID_{\delta_i},$
$\{(\phi_k,\psi_k)\}$ strongly converges to $(\phi,\psi)$ in
$W^{1,2}\times L^4$, so 
equivalently we need to prove that there are Dirac-harmonic spheres 
$(\sigma_i^l,\xi_i^l): \S^2\to
\S^n, i=1,2,\cdots,I; $ $l=1,2,\cdots,L_i$ such that
\begin{equation}
\label{c6} \Sigma^I_{i=1} {\rm lim }_{\delta_i\to 0}{\rm
lim}_{k\to\infty}E(\phi_k;D_{\delta_i})=\Sigma^I_{i=1}
\Sigma_{l=1}^{L_i}E(\sigma^l_i),
\end{equation}
\begin{equation}
\label{c7} \Sigma^I_{i=1} {\rm lim }_{\delta_i\to 0}{\rm
lim}_{k\to\infty}E(\psi_k;D_{\delta_i})=\Sigma^I_{i=1}
\Sigma_{l=1}^{L_i}E(\xi^l_i).
\end{equation}
In fact, we will prove for each $i=1,2,\cdots,I$ that
there exist Dirac-harmonic spheres 
$(\sigma_i^l,\xi_i^l): \S^2\to
\S^n, $ $l=1,2,\cdots,L_i$ such that

\begin{equation}
\label{c8}
 {\rm lim }_{\delta\to 0}{\rm
lim}_{k\to\infty}E(\phi_k;D_{\delta})=
\Sigma_{l=1}^{L}E(\sigma^l),
\end{equation}
\begin{equation}
\label{c9} {\rm lim }_{\delta\to 0}{\rm
lim}_{k\to\infty}E(\psi_k;D_{\delta})= \Sigma_{l=1}^{L}E(\xi^l),
\end{equation}
where, for simplicity, we have dropped the subscript $i$ of
$D_{\delta_i}, L_i$ etc. and denote $p_i$ by $p$. Certainly, in
each $D_{\delta}$ there is only one blow-up point $p$.
\par
%1
Let us first prove (\ref{c8}) and (\ref{c9}) for the case that
there is only one bubble at the blow-up point $p$. Then, we need to prove
that there exists Dirac-harmonic sphere 
$(\sigma^1,\xi^1): \S^2\to
\S^n, $ such that
\begin{equation}
\label{c8-1}
 {\rm lim }_{\delta\to 0}{\rm
lim}_{k\to\infty}E(\phi_k;D_{\delta})=
E(\sigma^1),
\end{equation}
\begin{equation}
\label{c9-1} {\rm lim }_{\delta\to 0}{\rm
lim}_{k\to\infty}E(\psi_k;D_{\delta})= E(\xi^1).
\end{equation}

\par
For each $(\phi_k,\psi_k)$, we choose $\lambda_k$ such that
$${\rm max}_{x\in D_\delta(p)} E(\phi_k,\psi_k;
D_{\lambda_k}(x))=\frac{\varepsilon_0}{2},$$ and then choose
$x_k\in D_\delta$ such that
 $$ E(\phi_k,\psi_k;
D_{\lambda_k}(x_k))=\frac{\varepsilon_0}{2}.
$$
We may assume that $\lambda_k \to 0$ and $x_k \to p$ as $k\to
\infty$. Let $$\widetilde{\phi}_k(x):=\phi_k(x_k+\lambda_kx),\quad
\widetilde{\psi}_k(x):=\lambda_k^{-\frac 1 2 }
\psi_k(x_k+\lambda_kx),$$ then
$$E(\widetilde{\phi}_k,\widetilde{\psi}_k;D_1(0))=E(\phi_k,\psi_k;D_{\lambda_k}(x_k))
=\frac{\varepsilon_0}{2}<\varepsilon_0,$$
$$E(\widetilde{\phi}_k,\widetilde{\psi}_k;D_K)=E(\phi_k,\psi_k;D_{\lambda_kK}(x_k))\leq
\Lambda,$$ from the $\varepsilon-$regularity theorem Theorem 3.2,
we have a subsequence, still denoted by
$(\widetilde{\phi}_k,\widetilde{\psi}_k)$ , that strongly
converges to some $(\widetilde{\phi},\widetilde{\psi})$ in
$W^{1,2}(D_R,N)\times L^4(\Sigma D_R\times \RR^K)$ for any $R\geq
1$. Thus, we get a nonconstant Dirac-harmonic map
$(\widetilde{\phi},\widetilde{\psi})$ on $\RR^2$, and by
stereographic projection, we obtain a nonconstant Dirac-harmonic
map $(\sigma^1,\xi^1)$ on $\S^2\setminus\{N\}$ with bounded
energy. By the regularity theorem, we have a nonconstant
Dirac-harmonic map $(\sigma^1,\xi^1)$ on the whole $\S^2,$ this is
the first bubble at the blow-up point $p$.
\par
Let $$A(\delta, R,k):=\{x\in \RR^2|\lambda_kR\leq
|x-x_k|\leq\delta\},$$ then to prove the assertion of (\ref{c8-1}) and (\ref{c9-1}) is equivalent to
proving 
\begin{equation}
\label{c12} {\rm lim}_{R\to\infty}{\rm lim}_{\delta\to 0}{\rm
lim}_{k\to \infty}E(\phi_k;A(\delta,R,k))=0,
\end{equation}
\begin{equation}
\label{c13} {\rm lim}_{R\to\infty}{\rm lim}_{\delta\to 0}{\rm
lim}_{k\to \infty}E(\psi_k;A(\delta,R,k))=0.
\end{equation}
Now let $(r_k,\theta_k)$ be the polar coordinate system centered
at $x_k$, $f: \RR\times \S^1\to \RR^2, \quad
f(t,\theta)=(e^{-t},\theta)$, where $\RR\times \S^1$ is given the
metric $g=dt^2+d\theta^2$, which is conformal to $ds^2$ on
$\RR^2$: $(f^{-1})^*g=e^{2t}ds^2.$ Denote $\Phi_k:=f^*\phi_k,$
$\Psi_k:=f^*\psi_k,$ then $E(\Phi_k,\Psi_k)\leq \Lambda$ since the
energy functional is conformally invariant. Denote $T_0:=|{\rm
log}\delta|,$ $P_T:=[T_0,T_0+T]\times\S^1$, $T>0$, then
$(\Phi_k,\Psi_k)\to (f^*\phi,f^*\psi):=(\Phi,\Psi)$ strongly on
$P_T$ for any $T>0$, hence $$E(\Phi_k,\Psi_k;P_T)\to
E(\phi,\psi;D_\delta\setminus D_{\delta e^{-T}}), \quad {\rm as}
\quad k\to \infty.$$ Given any small $\varepsilon >0$, since
$E(\phi,\psi)\leq \Lambda,$ we may choose $\delta >0$ small enough
such that $E(\phi,\psi;D_\delta)<\varepsilon/2,$ then for any
$T>0$, there is a $k(T)$ such that when $k\geq k(T)$,
\begin{equation}
\label{c14} E(\Phi_k,\Psi_k;P_T)<\frac{\varepsilon}{2}.
\end{equation}
Similarly, denote $T_k:=|{\rm log}\lambda_k R|,$
$Q_{T,k}:=[T_k-T,T_k]\times \S^1,$ then
\begin{equation}
\label{c15} E(\Phi_k,\Psi_k;Q_{T,k})<\frac{\varepsilon}{2}, \quad
k\geq k(T).
\end{equation}
Now we prove that there is a $K>0$ such that if $k\geq K$, then
\begin{equation}
\label{c16}
\int_{[t,t+1]\times\S^1}(|d\Phi_k|^2+|\Psi_k|^4)<\varepsilon,
\quad \forall t\in [T_0,T_k-1].
\end{equation}
Suppose this is false, then there exists a sequence of $\{t_k\}$
such that $t_k\to \infty$ as $k\to \infty,$ and
$$
\int_{[t_k,t_k+1]\times\S^1}(|d\Phi_k|^2+|\Psi_k|^4)\geq
\varepsilon, \quad \forall t\in [T_0,T_k-1],
$$
because of (\ref{c14}) and (\ref{c15}), we know that
$t_k-T_0,T_k-t_k\to \infty,$ by a translation from $t$ to $t-t_k$,
we get some $(\widetilde{\Phi}_k,\widetilde{\Psi}_k)$, and  for
all $k$ the following holds
\begin{equation}
\label{c17} \int_{[0,1]\times
\S^1}(|d\widetilde{\Phi}_k|^2+|\widetilde{\Psi}_k|^4)\geq
\varepsilon.
\end{equation}
By our assumption on $(\phi_k,\psi_k)$, we may assume that
$(\widetilde{\Phi}_k,\widetilde{\Psi}_k)$ weakly converges to
$(\Phi_\infty,\Psi_\infty)$ in $W_{loc}^{1,2}\times
L^4_{loc}(\RR\times\S^1)$. If this is strong convergence on
$[0,1]\times \S^1$, then we obtain a nonconstant Dirac-harmonic
map $(\Phi_\infty,\Psi_\infty)$ on the whole $\RR\times \S^1$, and
hence a nonconstant Dirac-harmonic map
$(\sigma_\infty,\xi_\infty)$ on $\S^2\slash\{N,S\}$ with finite
energy. Again by the regularity theorem, we have a nonconstant
Dirac-harmonic map $(\sigma_\infty,\xi_\infty)$ on $\S^2$, a
contradiction to the assumption that $L=1$. On the other hand, if
$(\widetilde{\Phi}_k,\widetilde{\Psi}_k)$ does not strongly
converge to $(\Phi_\infty,\Psi_\infty)$ on $[0,1]\times \S^1$,
then we may find some point $(t_0,\theta_0)\in[0,1]\times \S^1$ at
which $\{(\Phi_k,\Psi_k)\}$ blows up, in this case, we can still
get a second nonconstant Dirac-harmonic map
$(\sigma_\infty,\xi_\infty)$ on $\S^2$, again contradicting $L=1$.
Therefore, (\ref{c16}) holds true.
\par
Now from (\ref{c16}) and the $\varepsilon-$regularity theorem, we
have $$\|d\Phi_k\|_{L^\infty([t,t+1]\times\S^1)}\leq
C\sqrt{\varepsilon},$$ where $C>0$ is a constant independent of
$t$. Back to $\RR^2$, we have
\begin{equation}
\label{c18} |d\phi_k|(x)\leq \frac{C\sqrt{\varepsilon}}{r}, \quad
r=r(x)=|x|,\quad x\in A(\delta,R,k),
\end{equation}
from which we can conclude that
\begin{equation}
\label{c19} ||d\phi_k||_{L^{(2,\infty)}(A(\delta,R,k))}\leq
C\sqrt{\varepsilon},
\end{equation}
where $||\cdot ||_{L^{(2,\infty)}(A(\delta,R,k))}$ denotes the
norm in $L^{(2,\infty)}$ (cf. \cite{H1}).
\par
On the other hand, we  know that (for the target $\S^n$):
\begin{eqnarray*}
\Delta \phi^m&=&[(\phi^i_1\phi^m-\phi^i\phi^m_1)-\la
e_1\cdot\psi^i,\psi^m\ra_{\Sigma M}\phi^i_1]\nonumber
\\
&& +[(\phi^i_2\phi^m-\phi^i\phi^m_2)-\la
e_2\cdot\psi^i,\psi^m\ra_{\Sigma M}\phi^i_2]
\end{eqnarray*}
belongs to the Hardy space $\H^1$, so we have
\begin{equation}
\label{c20} ||d\phi_k||_{L^{(2,1)}(A(\delta,R,k))}\leq C,
\end{equation}
where $||\cdot||_{L^{(2,1)}(A(\delta,R,k))}$ denotes the norm in
the $L^{(2,1)}$ space. Therefore, by the duality of $L^{(2,1)} $
and $L^{(2,\infty)}$,  (\ref{c19}) and (\ref{c20}) we have
$$\int_{A(\delta,R,k)}|d\phi_k|^2\leq ||d\phi_k||^2_{L^{(2,\infty)}(A(\delta,R,k))}
||d\phi_k||^2_{L^{(2,1)}(A(\delta,R,k))}\leq C\varepsilon,$$ that
is, for $k,R$ large enough, $\delta$ small enough,
$$E(\phi_k;A(\delta,R,k))<C\varepsilon,$$
which proves (\ref{c12}).
\par
Now we turn to the proof of (\ref{c13}). Choose a cut-off function
$\eta$ on $D(x_k,2\delta)$ as follows:
$$\eta\in C^\infty_0(D_{2\delta}\setminus
D_{\lambda_kR/2});\qquad \eta\equiv 1 \quad {\rm in}\quad
D_{\delta}\setminus D_{\lambda_kR}$$
$$|\nabla\eta|\leq C/\delta \quad {\rm in}\quad D_{2\delta}\setminus
D_\delta;\qquad |\nabla\eta|\leq C/\lambda_kR \quad {\rm in}\quad
D_{\lambda_kR}\setminus D_{\lambda_kR/2},$$ where we denote
$D_\delta:=D(x_k,\delta)$ etc. for simplicity.
 Then from $\partial \hskip
-2.2mm \slash (\eta \psi_k)=$ $\eta \partial \hskip -2.2mm \slash
\psi_k+\nabla\eta\cdot\psi_k,$ we have
\begin{eqnarray*}
\|\eta\psi_k\|_{L^4}&\leq&C\|\eta \partial \hskip -2.2mm \slash
\psi_k+\nabla\eta\cdot\psi_k\|_{L^{\frac 4 3}} \\
&\leq&C\||\eta||d\phi_k||\psi_k|+|\nabla\eta||\psi_k|\|_{L^{\frac
4 3}}\\ &\leq&
C\|d\phi_k\|_{L^2(A(2\delta,R/2,k))}\|\psi_k\|_{L^4(A(2\delta,R/2,k))}
+C[\int_{A(2\delta,R/2,k)}(|\nabla\eta||\psi_k|)^{\frac 4
3}]^{\frac 3 4}\\ &\leq&
C\sqrt{\varepsilon}+[\int_{D_{2\delta}\setminus
D_\delta}(|\nabla\eta| |\psi_k|)^{\frac 4 3}]^{\frac 3
4}+[\int_{D_{\lambda_kR}\setminus
D_{\lambda_kR/2}}(|\nabla\eta||\psi_k|)^{\frac 4 3}]^{\frac 3 4},
\end{eqnarray*}
where in the last line we used (\ref{c12}). By the definition of
$\eta$ we have
\begin{eqnarray*}
\|\psi_k\|_{L^4(A(\delta,R,k))}&\leq& C\sqrt{\varepsilon}+ C
[\int_{D_{2\delta}\setminus D_\delta} |\psi_k|^4]^{\frac 1
4}+C[\int_{D_{\lambda_kR}\setminus
D_{\lambda_kR/2}}|\psi_k|^4]^{\frac 1 4}\\
&\leq&C\sqrt{\varepsilon}+C\varepsilon^{\frac 1
4}+C\varepsilon^{\frac 1 4},
%+\frac c \delta
%(\int_{D_{2\delta}-D_\delta}
%|\psi|^4)^{\frac 1 4}\delta^{\frac 3 2}+\frac{C}{\lambda_kR}
%(\int_{D_{\lambda_kR}-D_{\lambda_kR/2}}|\psi|^4)^{\frac 1 4}(\lambda_kR)^{\frac 3 2}\\
%&\leq& (\sqrt{\varepsilon}+\sqrt{\delta}+\sqrt{\lambda_kR}),
\end{eqnarray*}
where in the last step, we used (\ref{c16}).  This proves
(\ref{c13}). Therefore, the energy identities holds true for the
case $L=1$.
\par
For a fixed blow-up point $p$, the number $L$ of bubbles
$(\sigma,\xi)$ must be finite. This follows easily from the fact
that there is a number $C(\S^n)>0$ such that for all nonconstant
Dirac-harmonic maps $(\sigma,\xi): \S^2\to \S^n$, we have
$$E(\sigma,\xi;\S^2)\geq C(\S^n).$$ In fact, by Theorem 3.1, there
is a  constant $\varepsilon_0>0$ such that if
$E(\sigma,\xi;\S^2)<\varepsilon_0$, then $\sigma \equiv const. $
and $\psi$ satisfies the Dirac equation $\partial \hskip -2.2mm
\slash \psi=0$ on $\S^2$ which implies that $\psi\equiv 0$. We
therefore have $E(\sigma,\xi;\S^2)\geq \varepsilon_0$ for all
nonconstant Dirac-harmonic maps $(\sigma,\xi):\S^2\to \S^n$.
\par
The case of $L>1$ can be proved by induction on the number $L$, we
omit the details, as one may see the argument in \cite{DT}. \hfill
$\Box$

 \vskip12pt \noindent {\it Remark.} From the proof we see
that at each blow-up point $p_i$ $(i=1,2,\cdots,I)$, the
Dirac-harmonic maps $(\sigma^l_i,\xi^l_i): \S^2\to\S^n,
l=1,2,\cdots,L_i$ in Theorem 3.6 come from the blow-up process at
$p_i$.

\

\noindent{\it Acknowledgement.} The first author would like to thank the Max
  Planck Institute for Mathematics in the Sciences for support and
  good working conditions during the preparation of this paper. 
  We also would like to thank the referee for his or her careful reading.

\end{document}